\documentclass[11pt,a4paper, twoside]{article}
\usepackage[notref,notcite]{showkeys}

\usepackage{amsmath,amssymb,amsthm,amsfonts,mathrsfs,amscd,environ}
\usepackage{dsfont}
\usepackage{latexsym,enumerate,color,geometry,extarrows}
\usepackage{xcolor}
\usepackage{verbatim,fancyhdr,enumitem}
 \NewEnviron{ews}{%
\begin{equation}\begin{split}
  \BODY
\end{split}\end{equation}
}

\NewEnviron{ews*}{%
\begin{equation*}\begin{split}
  \BODY
\end{split}\end{equation*}
}

\def\beg{\begin}
\def\bequ{\begin{equation}}
\def\enqu{\end{equation}}
\def\bes{\begin{split}}
\def\ens{\end{split}}
\def\bews{\begin{ews}}
\def\beqn{\begin{eqnarray}}
\def\enqn{\end{eqnarray}}
\def\beq*{\begin{equation*}}
\def\enq*{\end{equation*}}
\def\bqn*{\begin{eqnarray*}}
\def\eqn*{\end{eqnarray*}}
\def\bary{\begin{array}}
\def\eary{\end{array}}
\def\bpma{\begin{pmatrix}}
\def\epma{\end{pmatrix}}
\def\bvma{\begin{Vmatrix}}
\def\evma{\end{Vmatrix}}

 \numberwithin{equation}{section}

\def\al{\alpha}
\def\be{\beta}
\def\ga{\gamma}
\def\de{\delta}
\def\ep{\epsilon}

\def\et{\eta}
\def\th{\theta}

\def\ka{\kappa}
\def\la{\lambda}
\def\rh{\rho}

\def\si{\sigma}
\def\ta{\tau}
\def\ph{\phi}

\def\ps{\psi}

\def\Si{\Sigma}
\def\Ph{\Phi}

\def\Om{\Omega}

\def\R{\mathbb R}
\def\P{\mathbb P}
\def\E{\mathbb E}

\def\N{\mathbb N}

\def\sF{\mathscr F}

\def\sB{\mathscr B}

\def\cW{\mathcal W}

\def\cT{\mathcal T}
\def\cM{\mathcal M}

\def\e{\operatorname{e}}

\def\d{\mathrm{d}}

\def\ff{\frac}
\def\ra{\rightarrow}
\def\nn{\nabla}
\def\pp{\partial}
\def\<{\langle}
\def\>{\rangle}
\def\sq{\sqrt}
\def\tld{\tilde}
\def\we{\wedge}
\def\1{\mathds{1}}

\def\trac{\mathrm{tr}}

\def\supp{\displaystyle{\mathrm{supp}}}

\def\hh{\hat }

\allowdisplaybreaks

\title{{\bf A Zvonkin's transformation for stochastic differential equations with singular drift and related applications}
}

\author{
{\bf Chenggui Yuan$^{a)}$~and~Shao-Qin Zhang$^{b)}$}\\
~\\
\footnotesize{$^{a)}$Department of Mathematics,   Swansea University, Bay Campus SA1 8EN, UK.}\\
\footnotesize{ Email: C.Yuan@swansea.ac.uk}\\
\footnotesize{$^{b)}$School of Statistics and Mathematics}\\
\footnotesize{Central University of Finance and Economics, Beijing 100081, China}\\
\footnotesize{Email: zhangsq@cufe.edu.cn}\\
}

\begin{document}

\maketitle

\begin{abstract}
In this paper, by establishing the $L^p$-$L^q$ estimate and Sobolev estimates for parabolic partial differential equations with a singular first order term  and a Lipschitz first order term, a new Zvonkin-type transformation is given for stochastic differential equations with singular  and   Lipschitz drifts. The associated Krylov's estimate is established. As applications, Harnack inequalities are established for stochastic equations with H\"older continuous diffusion coefficient and singular drift term without regularity assumption.    

\end{abstract}\noindent

AMS Subject Classification (2010): 60H10
\noindent

Keywords:  Zvonkin's transformation; singular diffusion processes; Krylov's estimate; Harnack inequality

\vskip 2cm

\section{Introduction}

In \cite{Zv}, a transformation that removes the drift of stochastic differential equation (in short SDE) was introduced by Zvonkin. This  transformation of the phase space  together with  Krylov's estimate (see \cite{Kry-86}) gives  a powerful tool  in studying  SDEs with irregular coefficients.  For instance, in \cite{Ve} the author first proved the existence and uniqueness of strong solutions to SDEs with bounded measurable drift; \cite{GM} proved the uniqueness of strong solution to SDEs with locally Lipschitz and strong elliptic diffusion coefficients and integrable drifts; \cite{ZX05} extended results to equations with  local integral  drifts which has linear growth and Sobolev diffusion coefficients.  Recently, \cite{KR} obtained the existence and uniqueness of strong solutions to SDEs with additive noise and  time dependent drifts satisfying the $L^p$-$L^q$ integration condition, see \eqref{def-lp-lq} for instance. Krylov and R\"ockner's results were extended by \cite{ZhangX11} to the case of multiplicative noise, and stochastic homeomorphism flow property of singular SDEs were studied therein. For more properties of singular SDEs investigated by using Zvonkin's transformation and Krylov's estimate, see \cite{HW,LLW,XZ16,XZ,ZhangX16} and reference therein. 

We consider the following equation 
\beg{align}\label{equ-main-0}
\d X_t= b(t,X_t)\d t+b_0(t,X_t)\d t+\si(t,X_t)\d W_t,
\end{align}
where $b(t,\cdot):\R^d\ra \R^d$ is Lipschitz uniformly w.r.t $t\geq 0$, $b_0:[0,\infty)\times\R^d\ra\R^d$ is singular term satisfying the $L^p$-$L^q$ condition as in  \cite{KR}: 
\beg{align}\label{def-lp-lq}
\int_0^T\left(\int_{\R^d} |b_0|^p(t,x)\d x\right)^{\ff q p}\d t<\infty,~T\in [0,\infty),
\end{align}
with $p,q\in (1,\infty)$ and $\ff d p +\ff 2 q<1$, $\si:[0,\infty)\times\R^d\ra\R^d\otimes\R^d$ is non-degenerate, $\{ W_t\}_{t\geq 0}$ is a Brownian motion w.r.t.  a probability space with filtration $(\Om,\sF,\{\sF_t\}_{t\geq 0},\P)$.  In this paper, we shall give a new Zvonkin-type transformation $\Ph_t(\cdot):\R^d\ra\R^d$  by solving a PDE associated with \eqref{equ-main-0}. Precisely, let $\ph(t,x)=(\ph^1(t,x),\cdots,\ph^d(t,x))$ satisfy
\beg{align}\label{n-ph0}
\pp_t \ph^{i}+\ff 1 2\trac(\si\si^* \nn^2\ph^{i})+\<b+b_0,\nn \ph^i\>= -b_0^i+\la \ph^i,~i=1,\cdots, d.
\end{align}
Then  $\Ph_t(x):=\ph(t,x)+x$ satisfies the following equation equivalently
\beg{align*}
\pp_t \Ph^{i}+\ff 1 2\trac(\si\si^* \nn^2\Ph^{i})+\<b+b_0,\nn \Ph^i\>= b^i+\la \ph^i,~i=1,\cdots, d.
\end{align*}
The equation \eqref{n-ph0}  is different from the parabolic equation considered  in \cite{XXZZ,XZ16,XZ,ZX05,ZhangX11,ZhangX16} and can not be covered by their studies since the coefficient $b$  can have linear growth.  In fact, the $L^p$-$L^q$ estimate established in \cite[(10.3)]{KR} or \cite[(3.2)]{XXZZ} fails for $\pp_t \ph$, see Theorem \ref{lem-lplq} and Remark \ref{rem-lplq} below. We solve \eqref{n-ph0} (in fact, a more general parabolic equation) in a weighted space, and more details on the well-posedness and a priori estimates are available in Theorem \ref{lem-lplq} below.  

We prove that $\Ph_t$ is a homeomorphism on $\R^d$  by choosing $\la$ large enough, see Theorem \ref{lem-lplq} and \eqref{Ph} below.   Let $\Ph_t^{-1}$ is the inverse of $\Ph_t$. By It\^o's formula (see Lemma \ref{ITO} for a proof), we have
 \beg{align*}
 \d \Ph_t(X_t) & = b(t,X_t)\d t+\la\ph(t,X_t)\d t+\mbox{martingale part}\\
 & = b(t,\Ph_t^{-1}(\Ph_t(X_t)))\d t+\la\ph_t(\Ph_t^{-1}(\Ph_t(X_t)))\d t+\mbox{martingale part}.
 \end{align*}
Then $b(t,\Ph_t^{-1}(\cdot))$, as a drift term of a SDE for $\Ph_t(X_t)$, remains to be Lipschitz. Moreover,   if $b$ is monotone  in addition:
\beg{align}\label{mono}
\<b(t,x)-b(t,y),x-y\>\leq K|x-y|^2,~x,y\in\R^d,
\end{align}
then  $b(t,\Ph_t^{-1}(\cdot))$ also satisfies \eqref{mono} with another constant of the same sign with $K$.  However, by applying Zvonkin's transformation used in \cite{XZ16,XZ,ZX05,ZhangX11,ZhangX16} to \eqref{equ-main-0}, one  gets a SDE with a locally Lipschitz drift term.  This property allows us to establish Harnack inequalities for \eqref{equ-main-0}. Harnack inequalities for SDEs with singular drifts have been  investigated in \cite{Huang,HZ,LLW,Shao}.  In \cite{LLW}, only log-Harnack inequality is established for SDEs with the drift satisfying the $L^p$-$L^q$ condition.  \cite{Shao} obtains Harnack inequalities with an  extra constant.  In \cite{Huang}, the author imposes extra regularities on  space variable, which turns out to requiring that the drift term should be H\"older continuous, see Remark \ref{rem-Har}. However, the drift term of the SDEs concerned about in \cite{Huang,LLW,Shao} can not include a Lipschitz drift. \cite{WangJDE} introduced a transformation for SDEs with Dini-continuous drift that retains the linear drift which automatically is Lipschitzian. Following this transformation, \cite{HZ} obtained Harnack inequalities for stochastic functional partial differential equations with Dini-continuous drift.  We establish Harnack inequality with power for \eqref{equ-main-0} under $L^p$-$L^q$ integral condition with $\ff d p +\ff 2 q<1$ and  the diffusion coefficient that can be H\"older continuous with order in $[\ff 1 2, 1]$.  Moreover, if $\ff d p+\ff 2 q<\ff 1 2$ and the H\"older continuity order of the  diffusion coefficient  is in $(\ff 1 2,1]$, then the Harnack inequality without  extra constant is established. We use a coupling modified from \cite{W11,WB} so that the diffusion coefficient can be H\"older continuity with index in $(\ff 1 2,1]$. This is new even in the case that the drift is regular.  

This paper is structured as follows. In Section 2, we investigate well-posedness and a priori estimates of a general parabolic equation which covers \eqref{n-ph0}.  Then Krylov's estimates for the solution of \eqref{equ-main-0} will be given in Section 3. In Section 4, we study Harnack inequality for the associated transition semigroup generated by \eqref{equ-main-0}.

Throughout this paper, we denote by $\|\cdot\|$ the operator norm of matrixes. For a (real, vector or matrix value)  function on $ [0,T]\times\R^d $, we denote  
\beg{align*}
\|f\|_{T,\infty} & =\sup_{t\in [0,T],x\in\R^d}\|f(t,x)\|,\\
\|f(t,\cdot)\|_\infty&=\sup_{x\in\R^d}\|f(t,x)\|.
\end{align*} 
Let $\{e_j\}_{j=1}^d$ be the ONB of $\R^d$. For any $A\in\R^d \otimes \R^d$, we denote $A_j^i=\<Ae_j,e_i\>$. For any $g\in C^1(\R^d)$, we denote by $\nn g(x)$ the gradient of $g$ at $x$ with 
$$(\nn g)^j(x):=\<\nn g(x),e_j\>:=\nn_{e_j} g(x).$$
For any $g\in C^1(\R^d,\R^d)$, we denote by $\nn g(x)\in \R^d\otimes \R^d$ the gradient of $g$ with
\beg{align*}
g^j(x):=\<g(x),e_j\>,\qquad (\nn g)^{j}_i(x) := \<\nn g(x)e_i,e_j\>=\nn_{e_i}g^j(x).
\end{align*}  
Particularly, for any $g\in C^2(\R^d)$, we denote by $\nn^2 g(x)\in \R^d\otimes\R^d$ the Hessen matrix of $g$ at $x$ with
$$(\nn^2 g)^i_j(x):=\<\nn^2 g(x)e_j,e_i\>:=\nn_{e_j}\nn_{e_i}g(x).$$

\section{$L^p$-$L^q$ estimates for parabolic equations}

We first study the $L^p$-$L^q$ estimates of the following parabolic equation
\beg{align}\label{equ-pa-1}
\pp_t u &+\trac{\left(a \nn^2 u \right)} + (b_1+b_2+b_0)\cdot\nn u + c u=\la u +f, 
\end{align} 
where the derivatives of $u$ are understood in the weak sense, and $a: [0,T]\times\R^d\ra \R^d\otimes\R^d$, $b_1,b_2,b_0: [0,T]\times\R^d\ra \R^d$ and $c,f: [0,T]\times\R^d\ra\R$ are measurable.  We assume that $a,b_1,b_2,c$ satisfy the following hypothesis.
\beg{description}[align=left, noitemsep]
\item [(H1)] $a$ is uniformly continuous in $x$ uniformly w.r.t. $t$, i.e. for any $T>0$ and $\ep>0$, there exists $\de>0$ such that  for any $x,y\in\R^d$ with $|x-y|<\de$
$$\sup_{t\in [0,T]}||a(t,x)-a(t,y)||<\ep.$$
For each $T>0$, there exist positive constants $\ka_1,\ka_2$ with $\ka_1\leq \ka_2$ such that
$$\ka_1|v|^2\leq \<a(t,x)v,v\>\leq \ka_2|v|^2,~(t,x)\in [0,T]\times\R^d,v\in\R^d.$$

\item [(H2)] For every $T>0$ and $t\in [0,T]$,  $b_1(t,\cdot)$ is Lipschitz continuous with Lipschitz constant $\|\nn b_1(t,\cdot)\|_\infty$, and   
$$\sup_{t\in [0,T]}\left(|b_1|(t,0)+\|\nn b_1(t,\cdot)\|_\infty\right)<\infty.$$

\item [(H3)] $b_2$ and $c$ are bounded on $[0,T]\times\R^d$
  
\end{description}

The condition {\bf (H2)} implies that $b_1$ has linear growth: there exists $K_0>0$ such that
\beg{align}\label{lin-gro}
\sup_{t\in [0,T]}|b_1(t,x)|\leq K_0(1+|x|),~x\in\R^d.
\end{align}
\beg{rem}\label{rem-0}
Let $b_1: [0,T]\times \R^d\ra \R^d$ satisfy $\|\nn b_1 \|_{T,\infty}<\infty$. 
Then there exist $\tld b_1(t,\cdot)\in C_b^2(\R^d)$ and bounded $\tld b_2$  such that $b_1=\tld b_1+\tld b_2$. In fact, for any nonnegative  $\et\in C^\infty_c(\R^d)$ with $\int_{\R^d}\et=1$, we set 
$$\tld b_1(t,x)=(b_1(t,\cdot)*\et)(x),\qquad \tld b_2(t,x)=b_1(t,x)-\tld b_1(t,x).$$
Then it is clear that $\tld b_1(t,\cdot)\in C^2(\R^d)$ and
\beg{align*}
 \|\nn \tld b_1\|_{T,\infty} +  \|\nn^2 \tld b_1 \|_{T,\infty} & \leq  \|\nn b_1 \|_{T,\infty},\\
 \|\tld b_2\|_{T,\infty}&\leq  C \|\nn b_1\|_{T,\infty}.
\end{align*}
\end{rem}
Due to this remark, we can use the following {\bf (H2')} instead of {\bf (H2)} under {\bf (H3)}:
\beg{description}
\item [(H2')] For every $T>0$ and $t\in [0,T]$,  $b_1(t,\cdot)\in C^2_b(\R^d)$ and 
$$\sup_{t\in [0,T]}\left(|b_1|(t,0)+\|\nn b_1(t,\cdot)\|_\infty+\|\nn^2 b_1(t,\cdot)\|_\infty\right)<\infty.$$
\end{description}

Before our main results, we introduce some function spaces which will be used through out this paper. Given $0\leq t\leq T$. Let  $w$ be a positive function, and let $L^p_w(\R^d) =L^p(\R^d, w(x)\d x)$,  $L^{p,w}_q(t,T)=L^q([t,T],L_w^p(\R^d))$, $L^{p}_q(t,T)=L^q([t,T],L^p(\R^d))$, and denote by $\|\cdot\|_{L^{p,w}}$, $\|\cdot\|_{L^{p,w}_q(t,T)}$,  $\|\cdot\|_{L^p_q(t,T)}$ the norms on these spaces respectively. We  define the Sobolev space $\left(W_{1,q}([t,T], L^p_w(\R^d)),\|\cdot\|_{W_{1,q}^{p,w}(t,T)}\right)$  as follows
\beg{align*}
W_{1,q}([t,T], L^p_w(\R^d))&=\left\{f\in L^{p}_q(t,T)~|~\pp_t f\in L^{p,w}_q(t,T)\right\},\\
\|f\|_{W_{1,q}^{p,w}(t,T)}&=\|\pp_t f\|_{L^{p,w}_q(t,T)}+\|f\|_{L^p_q(t,T)}.
\end{align*}
We denote by $\|\cdot\|_{W^{2,p}_q(t,T)}$ the norm of $ L^q([t,T],W^{2,p}(\R^d))$ where $W^{2,p}(\R^d)$ is the second-order Sobolev space. Let 
\beg{align*}
\cW^{2,p,w}_{1,q}(t,T)&=W_{1,q}([t,T], L^p_w(\R^d))\bigcap L^q([t,T],W^{2,p}(\R^d)),\\
\|f\|_{W^{2,p,w}_{1,q}(t,T)}&=\|f\|_{W_{1,q}^{p,w}(t,T)}+\|f\|_{W^{2,p}_q(t,T)}. 
\end{align*}
We  denote by $\cW^{2,p}_{1,q}(t,T)$ the case that $w\equiv 1$. Let $C^\al(\R^d)$ be the H\"older space on $\R^d$. We denote $C^{\al}_q(t,T)=L^{q}([t,T],C^{\al}(\R^d))$ with the norm
$$\|f\|_{C^\al_q(t,T)}=\left(\int_t^T \|f(t,\cdot)\|_{C^{\al}(\R^d)}^q\d t\right)^{\ff 1 q},~f\in C^\al_q(t,T).$$
If $t=0$, then $L^{p,w}_q(0,T)$, $L^p_q(0,T)$ e.t.c.  will be denoted by $L^{p,w}_q(T)$, $L^p_q(T)$ e.t.c.

We prove that \eqref{equ-pa-1} has a unique solution in the weighted space $\cW^{2,p,w}_{1,q}(t,T)$ with suitable weight $w$. In the following theorem, we denote $(\la-\la_0)^-= [ -(\la-\la_0)]\vee 0$.

\beg{thm}\label{lem-lplq}
Let $p,q\in (1,\infty)$ and $p_1\in [p,+\infty]$ with $\ff d {p_1}+\ff 2 q<1$.  Assume that $ b_0\in L^{p_1}_q(T)$, $  f\in L^{p}_q(T)$ and $u_0\in B^{2-2/q}_{p,q}(\R^d)$, and that {\bf (H1)}-{\bf (H3)} hold. Then\\
(1) \eqref{equ-pa-1} has a unique solution in $\cW_{1,q}^{2,p,w} (T)$ with $w(x)=\left(1+|x|^2\right)^{-\ff p 2} $.  Moreover, there exist  a constant $\la_0>0$ and  a positive constant  $C_1$ which depends on $\ka_1,\ka_2$, $p,q,d,T$, $\|\nn b_1\|_\infty$, $\| b_2\|_\infty,\| c\|_\infty$ and $(\la-\la_0)^-$ such that 
\beg{align}\label{ineq-vv}
&(\la\vee \la_0)\|u\|_{L^p_q(t,T)}+\|(\pp_t  +b_1\cdot\nn )u\|_{L^p_q(t,T)}+\|u\|_{W^{2,p}_q(t,T)}\nonumber\\
&\qquad\leq C_1\e^{C_1 \| b_0\|_{L^{p_1}_q(t,T)}^q}\left(\|  f\|_{L^p_q(t,T)}+\| u_T\|_{B^{2-2/q}_{p,q}}\right),~0\leq t\leq T.
\end{align}
(2) For any $\al\in [0,2)$, $p_2\in [p,+\infty]$, and $q_2\in [q,+\infty]$ such that 
$$\be_0:=\ff 1 2\left(2-\al+\ff 2 {q_2}+\ff d {p_2}-\ff 2 q-\ff d p\right)>0,$$
then there exists a constant $C_2>0$ depending on $\ka_1,\ka_2$, $p_2,q_2$,$p,q,d,T$, $\|\nn b_1\|_\infty$, $\| b_2\|_\infty,\| c\|_\infty$ and $(\la-\la_0)^-$ such that for any $0<\be<\be_0$
\beg{align}\label{add-la-be}
&\|u\|_{C^{\al}_{q_2}(t,T)}\1_{\{  p_2= +\infty \}}+\|u\|_{W^{\al,p_2}_{q_2}(t,T)}\1_{\{  p_2 \in  [p,+\infty) \}}\nonumber\\
&\qquad\qquad \leq \ff {C_2} {(\la\vee \la_0)^\be}\left(\|  f\|_{L^p_q(t,T)}+\| u_T\|_{B^{2-2/q}_{p,q}}\right),~0\leq t\leq T.
\end{align}
Moreover, if $q_2=q$ and $p_2\in [p,+\infty)$, then \eqref{add-la-be} holds with $\be=\be_0$.
\end{thm}

\beg{rem}\label{rem-lplq}
The \eqref{equ-pa-1} with $u(T,\cdot)=0$ is said to have $L^q$-Maximal regularity if 
$$\| \pp_t   u\|_{L^p_q(t,T)}+\|u\|_{W^{2,p}_q(t,T)}\leq C \|  f\|_{L^p_q(t,T)}.$$
Maximal $L^q$-regularity for evolution equations with time independent operators $A$ implies that $A$ generates an analytic semigroup, see \cite[Proposition 2.2]{Mo}. However, the generator of O-U semigroup can not generate an analytic semigroup in $L^p(\R^d)$, see \cite{LV}. Note that $b_1$ can has linear growth, which yields that $b_1(t,\cdot) \notin L^p(\R^d)$ and the elliptic operator in \eqref{equ-pa-1} covers the generator of O-U semigroup. Then one can not expect to derive the maximal regularity for \eqref{equ-pa-1} in $\cW^{2,p}_{1,q}(T)$. Combining this with \eqref{ineq-vv}, it is clear that $\|\pp_tu\|_{L^{p}_q(t,T)}$ can not be controlled by the right hand side of \eqref{ineq-vv}.  Hence, we combine $\pp_t u$ and $b_1\cdot \nn u$ together.

\end{rem}

\subsection{Proofs of Theorem \ref{lem-lplq}}
We first investigate \eqref{equ-pa-1} with $b_1\equiv 0$. 

\beg{lem}\label{lem-v}
Let $p,q\in (1,\infty)$ and $p_1\in [p,+\infty]$ with $\ff d {p_1}+\ff 2 q<1$.  Assume that $b_1\equiv 0$, $ b_0\in L^{p_1}_q(T)$, $  f\in L^{p}_q(T)$ and $u_0\in B^{2-2/q}_{p,q}(\R^d)$, and that {\bf (H1)} and {\bf (H3)} hold. Then all the assertions in Theorem \ref{lem-lplq} hold.
\end{lem}
\beg{rem}
For $b_0=0, b_1=0$, it has been proved in \cite[Theorem 1.2 and Theorem 5.4]{GV} that \eqref{equ-pa-1} has a unique solution for any $p,q\in(1,\infty)$ and \eqref{ineq-vv} holds. Then by the continuity method, to prove this lemma, it suffices to show \eqref{ineq-vv} assuming that the solution already exists.
\end{rem}
\beg{proof}
(1) Let $u\in \cW^{2,p}_{1,q}(T)$. We first consider a prior estimate of $ b_0\cdot\nn u$:
\beg{align}\label{add-bnnu}
\|  b_0(t,\cdot)\cdot\nn u(t,\cdot)\|_{L^p}\leq C\|  b_0(t,\cdot)\|_{L^{p_1}}\| u (t,\cdot)\|_{B^{2-\ff 2 q}_{p,q}}.
\end{align}
By \cite[III Theorem 4.10.2]{Am} and  the interpolation
$$(L^p(\R^d),W^{2,p}(\R^d))_{1-\ff 1 q,q}=B_{p,q}^{2-\ff 2 q}(\R^d),$$  
we have that $u\in C([0,T], B^{2-\ff 2 q}_{p,q}(\R^d))$.
For $p_1>p$. Since $2-\ff 2 q>1+\ff d {p_1}$, we have by \cite[2.3.3/(9)]{Tr78}  that  $B^{2-\ff 2 q}_{p,q}(\R^d)$ continuously embeds into $W^{1+\ff d {p_1}}(\R^d)$. Thus 
$$\|\nn u(t,\cdot)\|_{L^{\ff {pp_1} {p_1-p}}}\leq \| u(t,\cdot)\|_{W^{1,\ff {pp_1} {p_1-p}}}\leq C\|u(t,\cdot)\|_{W^{1+\ff d {p_1},p}}\leq C\|u(t,\cdot)\|_{B^{2-\ff 2 q}_{p,q}}.$$
For $p_1=p$.  Since $2-\ff 2 q>1+\ff d p_1=1+\ff d p$, $B^{2-\ff 2 q}_{p,q}(\R^d) $ continuously embeds into the H\"older space $C^{1+\ep}(\R^d)$ with any $0<\ep<1-\ff 2 q-\ff d p$.  The embedding theorem yields that for any $\de>0$ 
$$\| \nn u(t,\cdot)\|_{L^{\infty}}\leq C\| u(t,\cdot)\|_{B^{2-\ff 2 q}_{p,q}}.$$
Setting $\ff {pp_1} {p_1-p}:=p$ if $p_1=+\infty$. Then,  for $p_1\geq p$,  by  H\"older's inequality,
\beg{align*}
\| b_0(t,\cdot)\cdot\nn u(t,\cdot)\|_{L^p}&\leq \|b_0(t,\cdot)\|_{L^{p_1}}\|\nn u(t,\cdot)\|_{L^{\ff {pp_1} {p_1-p}}}\\
&\leq C\| b_0(t,\cdot)\|_{L^{p_1}}\|  u(t,\cdot)\|_{B^{2-\ff 2 q}_{p,q}}.
\end{align*}

By \cite[III Theorem   4.10.2 and Lemma 4.10.1]{Am}, there exists $\tld C>0$ such that 
\beg{align}\label{add-B-cW}
\sup_{s\in [t,T]}\|  u(s,\cdot)\|_{B^{2-\ff 2 q}_{p,q}}\leq \tld C\|u\|_{\cW^{2,p}_{1,q}(t,T)},~0\leq t\leq T.
\end{align}
By \eqref{add-bnnu} and $u\in C([0,T], B^{2-\ff 2 q}_{p,q}(\R^d))$, we have that
\beg{align*}
\| b_0 \nn u\|_{L^p_q(t,T)}&\leq C\|b_0\|_{L^{p_1}_q(t,T)}\sup_{s\in [t,T]}\|  u(s,\cdot)\|_{B^{2-\ff 2 q}_{p,q}}\\
&\leq C\|b_0\|_{L^{p_1}_q(t,T)}\|u\|_{\cW^{2,p}_{1,q}(t,T)}<\infty,~0\leq t\leq T.
\end{align*}
Hence, $(\pp_t+L_t+c-\la)$ is a bounded operator from $\cW^{2,p}_{1,q}(T)$ to $L^p_q(T)$.

Next, we  prove \eqref{ineq-vv}. By \cite[Theorem 1.2 and Theorem 5.4]{GV}, there exists a constant $\la_0>0$ and a positive constant $C_0$ which depends on $p,d,q, T$, $\| b_2\|_\infty$, $\| c\|_\infty$ and $(\la-\la_0)^-$ such that for any $0\leq t\leq T$
\beg{align}\label{ineq-vf0}
(\la\vee \la_0)\|u\|_{L^p_q(t,T)}&+\|\pp_t u\|_{L^p_q(t,T)}+\|u\|_{W^{2,p}_q(t,T)}\nonumber\\
&\leq C_0\left(\| f+ b_0\cdot\nn u\|_{L^p_q(t,T)}+\|u_T\|_{B^{2-2/q}_{p,q}}\right).
\end{align}
Let 
$$I(t)=\left((\la\vee \la_0)\|u\|_{L^p_q(t,T)}+\|u\|_{\cW^{2,p}_{1,q}(t,T)}\right)^q.$$ 
Then it follows from \eqref{add-B-cW} and \eqref{add-bnnu} that
\beg{align*}
\| b_0\cdot\nn u\|_{L^p_q(t,T)}^q&=\int_t^T \| b_0(s,\cdot)\cdot \nn u(s,\cdot)\|_{L^p}^q\d s\\
&\leq C\int_t^T \|  u(s,\cdot)\|_{B^{2-2/q}_{p,q}}^q\| b_0(s,\cdot) \|_{L^{p_1}}^q\d s\\
& \leq C\tld C\int_t^T I(s)\| b_0(s,\cdot)\|_{L^{p_1}}^q\d s
\end{align*}
Putting this into \eqref{ineq-vf0}, we have
\beg{align*}
I(t)\leq C_1\left(\|  f\|_{L^p_q(t,T)}+\|u_T\|_{B^{2-2/q}_{p,q}}^q\right)+C_1\int_t^T I(s)\| b_0(s,\cdot) \|_{L^{p_1}}^q\d s
\end{align*}
for some $C_1$ depending on $p,d,q,T$, $\| b_2\|_\infty$, $\| c\|_\infty$, $(\la-\la_0)^-$. Therefore, \eqref{ineq-vv} follows by Gronwall's inequality. 

(2) We use interpolation theorems to prove this assertion.   By \eqref{ineq-vv}, 
\beg{align*}
\|u\|_{L^p_q(t,T)}\leq \ff {C_1\e^{C_1 \| b_0\|_{L^{p_1}_q(t,T)}^q}} {\la\vee \la_0}\left(\|  f\|_{L^p_q(t,T)}+\| u_0\|_{B^{2-2/q}_{p,q}}\right).
\end{align*}
Denote by 
$$D_1=C_1\e^{C_1 \| b_0\|_{L^{p_1}_q(t,T)}^q},~A=\|  f\|_{L^p_q(t,T)}+\| u_0\|_{B^{2-2/q}_{p,q}}.$$
Then
\beg{align*}
\|u\|_{L^p_q(t,T)}\leq \ff {D_1A} {(\la\vee \la_0)},\qquad \|u\|_{W^{2,p}_q(t,T)}\leq  D_1A  . 
\end{align*}
By the Gagliardo-Nirenberg inequalities (e.g. \cite[Theorem 1]{BM}), for any $\al\in (0,2)$, we have that
\beg{align}\label{add-Wapq}
\|u\|_{W^{\al, p}_q(t,T)}&\leq C\left(\int_t^T\|u(s,\cdot)\|_{W^{2, p}}^{q \al/2 }\|u(s,\cdot)\|_{L^p}^{q(1-\al/2)}\d s\right)^{\ff 1 q}\nonumber\\
&\leq C(D_1A)^{ \al/2}\left(\ff {D_1A} {(\la\vee \la_0)}\right)^{1-\al/2}\nonumber\\
&= CD_1 A (\la\vee \la_0)^{-\ff 1 2(2-\al)},
\end{align}
which yields \eqref{add-la-be} with $p_2=p$, $q_2=q$. By the Sobolev embedding theorem, 
\beg{align}\label{Walp2q}
\|u\|_{W^{\al, p_2}_q(t,T)}\leq C\|u\|_{W^{\al-\ff d {p_2}+\ff d p, p}_q(t,T)}\leq \ff {CD_1 A} { (\la\vee \la_0)^{\ff  1 2 (2-\al+\ff d {p_2}-\ff d p)}},
\end{align}
which is  \eqref{add-la-be} with $p_2\in (p,+\infty)$, $q_2=q$.

For any $\th\in (0,1)$, by \cite[III Theorem 4.10.2]{Am} and 
$$(L^p(\R^d),W^{2(1-\th),p}(\R^d))_{1-\ff 1 q,q}=B^{2(1-\ff 1 q)(1-\th)}_{p,q}(\R^d),$$
 we have that
 \beg{align*}
\sup_{s\in[t,T]}\|u(s,\cdot)\|_{B^{2(1-\ff 1 q)(1-\th)}_{p,q}}\leq C\left(\|u\|_{W^{2(1-\th),p}_q(t,T)}+ \|\pp_t u\|_{L^p_q(t,T)}\right).
\end{align*}
Letting $u_a(t,x)=u(at,x)$, we have that for any $a>0$ 
\beg{align*}
\|u(s,\cdot)\|_{B^{2(1-\ff 1 q)(1-\th)}_{p,q}}&=\|u_a(s/a,\cdot)\|_{B^{2(1-\ff 1 q)(1-\th)}_{p,q}}\\
&\leq C\left( \|u_a\|_{W^{2(1-\th),p}_q(t,T)}+ \|\pp_t( u_a)\|_{L^p_q(t,T)}\right)\\
&\leq  C\left( \|u(a\cdot,\cdot)\|_{W^{2(1-\th),p}_q(t,T)}+a\|(\pp_t u)(a\cdot,\cdot)\|_{L^p_q(t,T)}\right)\\
&\leq  C\left(a^{-\ff 1 q}\|u\|_{W^{2(1-\th),p}_q(t,T)}+a^{1-\ff 1 q}\|\pp_t u\|_{L^p_q(t,T)}\right),~t\leq s\leq T,
\end{align*}
where we use \cite[III Lemma 4.10.1]{Am} in the last inequality. Then, by choosing the optimal $a$, we arrive at
\beg{align*}
\sup_{s\in[t,T]}\|u(s,\cdot)\|_{B^{2(1-\ff 1 q)(1-\th)}_{p,q}}&=  C\|u\|_{W^{2(1-\th),p}_q(t,T)}^{1-\ff 1 q} \|\pp_t u\|_{L^p_q(t,T)}^{\ff 1 q}.
\end{align*}
Then, taking into account \eqref{add-Wapq} and \eqref{ineq-vv},
\beg{align*}
\sup_{s\in[t,T]}\|u(s,\cdot)\|_{B^{2(1-\ff 1 q)(1-\th)}_{p,q}}\leq CD_1A(\la\vee \la_0)^{- (1-\ff 1 q)\th}
\end{align*}
Setting $\al=2(1-\ff 1 q)(1-\th)$, we have that
\beg{align}\label{add-Walp00}
\sup_{s\in[t,T]}\|u(s,\cdot)\|_{B^{\al}_{p,q}}\leq \ff {CD_1A} {(\la\vee \la_0)^{(1-\ff 1 q-\al/2)}} =\ff {CD_1A} {(\la\vee \la_0)^{ \ff 1 2(2-\ff 2 q-\al)}},~\al\in(0,2-\ff 2 q).
\end{align}
Since for any $p_2\in [p,+\infty)$, $W^{\al,p_2}(\R^d)\subset B^{\al-\ff d {p_2}+\ff d p+\ep}_{p,q}(\R^d)$. Then, for $\ep>0$ with $\al-\ff d {p_2}+\ff d p+\ep<2-\ff 2 q$,  it follows from \eqref{add-Walp00} that
\beg{align}\label{add-Walp}
\|u\|_{W^{\al,p_2}_\infty(t,T)}\leq C\sup_{s\in[t,T]} \|u(s,\cdot) \|_{B^{\al-\ff d {p_2}+\ff d p+\ep,p}_{p,q}}\leq \ff {CD_1A} {(\la\vee \la_0)^{ \ff 1 2(2-\ff 2 q-\al+\ff d {p_2}-\ff d p-\ep)}},
\end{align}
which implies \eqref{add-la-be} with $q_2=+\infty$ and $p_2\in [p,+\infty)$. 

Let $\ga\in (0,1)$, $p'\in [p,\infty)$, $\al_0\in (0,2-\ff 2 q)$ and $\al_1\in [0,2)$. Then for 
\beg{align*}
\al_2=\ga\al_0+(1-\ga)\al_1,~q_2=\ff {q} {1-\ga},~\ff 1 {p_2}=\ff {\ga} {p}+\ff {1-\ga} {p'},
\end{align*}
we have by the  the Gagliardo-Nirenberg inequality that
\beg{align*}
\|u(t,\cdot)\|_{W^{\al_2,p_2}}\leq C\|u(t,\cdot)\|_{W^{\al_0,p}}^{\ga}\|u(t,\cdot)\|_{W^{\al_1,p'}}^{1-\ga}.
\end{align*}
Combining this with \eqref{add-Walp} and \eqref{Walp2q},  for any $\ep>0$ with $\al_0+\ep<2-\ff 2 q$, we have that
\beg{align*}
\|u \|_{W^{\al_2,p_2}_{q_2}(t,T)}&\leq C\|u\|_{W^{\al_0,p}_\infty(t,T)}^{\ga}\|u \|_{W^{\al_1,p'}_q(t,T)}^{ 1-\ga}\\
&\leq  \ff {CD_1A} {(\la\vee \la_0)^{ \ff {\ga} 2(2-\ff 2 q-\al_0-\ep)+\ff {1-\ga} {2} (2-\al_1+\ff d {p'}-\ff d p)}}\\
&=\ff {CD_1A} {(\la\vee \la_0)^{ \ff {1} 2(2-\ff {2\ga} q-\al_2-\ga\ep+\ff {(1-\ga) d} {p'}-\ff {(1-\ga)d} p)}}\\
&=\ff {CD_1A} {(\la\vee \la_0)^{ \ff {1} 2(2-\ff {2} q+\ff 2 {q_2}-\al_2-\ga\ep+\ff d {p_2}-\ff {d} p)}}.
\end{align*}
Hence \eqref{add-la-be} holds with $q_2\in (q,+\infty)$, $p_2\in [p,+\infty)$. 

Since $B^{\al+\ff d p+\ep}_{p,q}(\R^d)\subset C^{\al}(\R^d)$  for any $\ep>0$, we have  
\beg{align}\label{add-Cal1}
\sup_{s\in [t,T]}\|u(s,\cdot)\|_{C^{\al}}\leq \ff {CD_1A} {(\la\vee \la_0)^{ \ff 1 2(2-\ff 2 q-\al-\ff d p-\ep)}},~0\leq \al<2-\ff 2 q-\ff d p-\ep.
\end{align}
Since $W^{\al+\ff d {p}+\ep,p}(\R^d)\subset C^{\al}(\R^d)$ for any $\ep>0$, we have that for any $q_2\in [q_2,+\infty)$ and $\al\geq 0$ such that 
$$2-\ff 2 q+\ff 2 {q_2}>\al+\ff d p+\ep,$$
there is
\beg{align}\label{add-Cal2}
 \|u\|_{C^{\al}_{q_2}(t,T)} \leq C\|u\|_{W^{\al+\ff d {p}+\ep,p}_{q_2}(t,T)}\leq \ff {CD_1A} {(\la\vee \la_0)^{ \ff 1 2(2-\ff 2 q+\ff 2 {q_2}-\al-\ff d p-\ep)}}.
\end{align}
\eqref{add-Cal1} and \eqref{add-Cal2} implies \eqref{add-la-be} since $\ep>0$ is arbitrary.

\end{proof}

Next, we give a transform to remove $b_1$ under the assumption {\bf (H2')}. This transformation has been used in \cite{CF,MPV} to investigate elliptic operators with unbounded coefficients. Consider the following ordinary differential system
\beg{align}\label{char-lin}
\ff {\d \ps} {\d t}(t,x)= b_1(t,\ps(t,x)),~\ps(T,x)=x.
\end{align}
For the solution of \eqref{char-lin}, we have the following lemma.
\beg{lem}\label{lem-ps}
Assume that {\bf (H2')} holds. Then, for any $t\in [0,T]$, $\ps(t,\cdot)\in C_b^2(\R^d)$ and $\ps(t,\cdot)$ is  a diffeomorphism  on $\R^d$. Denote by $\ps^{-1}(t,\cdot)$ with the inverse of $\ps(t,\cdot)$. Then $\ps^{-1}$ satisfies the following  ordinary differential system
\beg{align}\label{inv}
\ff {\d \ps^{-1}}  {\d t}(t,x)=-(\nn\ps)^{-1}(t,\ps^{-1}(t,x))b_1(t,x),\qquad \ps^{-1}(T,x)=x
\end{align}
and $\ps^{-1}(t,\cdot)\in C^2_b(\R^d)$. Moreover, 
\beg{align}\label{sup-ps}
 \|\nn\ps \|_{T,\infty}&+\|(\nn\ps)^{-1} \|_{T,\infty}+\|\nn \ps^{-1} \|_{T,\infty}\nonumber\\
 & +\sup_{1\leq l\leq d}\left(\|\nn^2\ps^l \|_{T,\infty}+\|\nn^2[(\ps^{-1})^l]\|_{T,\infty}\right)<\infty,
\end{align}
where $(\nn \ps)^{-1}(t,x)$ is the inverse  of the matix $\nn \ps(t,x)$ and satisfies
\beg{align}\label{inv-inv}
(\nn\ps)^{-1}(t,x)=(\nn\ps^{-1})(t,\ps(t,x)),~(t,x)\in[0,T]\times\R^d,
\end{align}
and the upperbound of \eqref{sup-ps} only depends on $T,\|\nn b_1\|_{T,\infty},\|\nn^2 b_1\|_{T,\infty},d$. 
\end{lem}
The proof of this lemma is direct, and we omit it. 

For any $u\in C^{1,2}([0,T]\times \R^d)$,  by a direct calculus, we have that
\beg{align*} 
&(\nn u)(t,\ps(t,x))=\left\{[(\nn\ps)^{-1}]^*\nn v\right\}(t,x),\\
&(\nn^2u)(t,\ps(t,x))=\left\{[\left(\nn\ps\right)^{-1}]^*\nn^2 v\left[(\nn\ps)^{-1}\right[ \right\}(t,x)\nonumber\\
&\qquad-\sum_{j=1}^d\left\{[\left(\nn\ps\right)^{-1}]^*(\nn^2\ps^j) (\nn\ps)^{-1}\right\}\left([(\nn\ps)^{-1}]^*\nn v\right)^j(t,x). 
\end{align*} 
Then, for $v(t,x):=u(t,\ps(t,x))$
\beg{align*}
&\qquad~~~~\pp_t v(t,x)= \pp_t u(t,\ps(t,x))+(b_1\cdot\nn u)(t,\ps(t,x)),\\
&(L_t u)(t,\ps(t,x)) =\trac\left(a(t,\ps)[(\nn\ps)^{-1}]^*\nn^2 v(\nn\ps)^{-1}\right)(t,x) \\
&\qquad -\left\{\sum_{j=1}^d \trac\left(a(t,\ps)[\left(\nn\ps\right)^{-1}]^*(\nn^2\ps^j) (\nn\ps)^{-1}\right)\left([(\nn\ps)^{-1}]^*\nn v\right)^j\right\}(t,x)\\
&\qquad + \left\{\left[(\nn \ps)^{-1}(b_0+b_1+b_2)(t,\ps)\right]\cdot \nn v\right\}(t,x).
\end{align*}
Let $\bar L_t$ be a differential operator defined as follows
$$\bar L_t g=\trac(\bar a \nn^2 g) +(\bar b_2+\bar b_0)\cdot\nn g,~g\in C^2([0,T]\times \R^d),$$
where 
\beg{align*}
\bar a(t,x) & =\left\{(\nn\ps)^{-1}a(t,\ps)[(\nn\ps)^{-1}]^*\right\}(t,x),\\
\bar b_2(t,x)&=(\nn\ps)^{-1}(t,x)b_2(t,\ps(t,x))\\
 &-\sum_{j=1}^d\left\{\trac\left(a(t,\ps)[\left(\nn\ps\right)^{-1}]^*(\nn^2\ps^j) (\nn\ps)^{-1}\right) [(\nn\ps)^{-1}]_j \right\}(t,x)\\
 \bar b_0(t,x)&= (\nn\ps)^{-1} (t,x)b_0(t,\ps(t,x)).
\end{align*}
Then, by setting 
$$ \bar f(t,x)=f(t,\ps(t,x)),\qquad \bar c(t,x)=c(t,\ps(t,x)),$$ 
we have  that
\beg{align}\label{equ-pa-2}
[(\pp_t+\bar L_t&+\bar c-\la) v](t,x) -\bar f(t,x)\nonumber \\
&=[ (\pp_t+L_t+c-\la)u](t,\ps(t,x))- f(t,\ps(t,x)),\\
v(T,\cdot) & =u_T(\cdot).
\end{align}
It follows from Lemma \ref{lem-ps} and {\bf (H3)} that $\bar b_2, \bar c$ are bounded.  If $b_0\in L^p_q(T)$, then $\bar b_0\in L^p_q(T)$. Moreover, since {\bf (H1)}, \eqref{sup-ps} and \eqref{inv-inv}, it is clear that $\bar a$ is  uniformly continuous and uniformly elliptic with some $\bar\ka_1$, $\bar \ka_2$ such that
$$\bar \ka_1|w|^2\leq \<\bar a(t,x)w,w\>\leq \bar \ka_2|w|^2,~(t,x)\in [0,T]\times\R^d, w\in\R^d.$$
For any $g\in\sB([0,T]\times \R^d)$, we define a mapping $J$ as follows 
$$(Jg)(t,x)=g(t,\ps^{-1}(t,x)),~(t,x)\in [0,T]\times\R^d.$$
Hence,  due to  \eqref{equ-pa-2}, we can investigate \eqref{equ-pa-1} by applying Lemma \ref{lem-v} to 
\beg{align}\label{equ-pa-v}
(\pp_t +\bar L_t+\bar c-\la) v=\bar f,~v(T,\cdot)=u(T,\cdot),
\end{align}
and showing that $u=Jv$ satisfies  \eqref{equ-pa-1} and all the  assertions of Theorem \ref{lem-lplq}.

Let 
\beg{align*}
\cW^{2,p}_{1,q,b_1}(T)=\left\{g\in\cW^{2,p,w}_{1,q}(T)~\Big|~(\pp_t+b_1\cdot\nn )g \in L^p_q(T) \right\}
\end{align*}
equipped the norm  
$$\|\cdot\|_{\cW^{2,p,w }_{1,q,b_1}(T)}:=\|(\pp_t+b_1\cdot\nn )\cdot\|_{L^p_q(T)}+\|\cdot\|_{\cW^{2,p,w}_{1,q}(T)}.$$ 
The following lemma shows that $J$ has nice properties.

\beg{lem}\label{lem-J}
Let $p,q\in (1,\infty)$ and  $w(x)=(1+|x|^2)^{-\ff p 2}$.\\
(1) Define $\{J_t\}_{[0,T]}$ as follows
$$(J_th)(x)=h(\ps^{-1}(t,x)),~h\in\sB(\R^d).$$
Then $\{J_t\}_{[0,T]}$ are uniformly bounded linear operators on $W^{j,p}(\R^d)$ with any $j\in \{0,1,2\}$. Consequently, $J$ is a homeomorphism on $W^{2,p}_q(T)$. Moreover, for any $g\in W^{2,p}_q(T)$ and a.e. $t\in [0,T]$, the following equality holds in $L^p(\R^d)$
\beg{align}\label{nn-u-v-1}
(\nn Jg)(t,\cdot)&=\left[J[(\nn\ps)^{-1}]^*\nn g)\right](t,\cdot),
\end{align}
(2) $J$ is a bounded operator from $\cW^{2,p }_{1,q}(T)$ to $ \cW^{2,p,w}_{1,q}(T)$, and
\beg{align}\label{pu-pv-bnu}
(\pp_t+b_1\cdot\nn)Jg&=J\pp_t g,~g\in\cW^{2,p }_{1,q}(T).
\end{align}
Moreover, $J$ is a homeomorphism from $\cW^{2,p }_{1,q}(T)$ to $ \cW^{2,p }_{1,q,b_1}(T)$, and there exist positive constants $c_1,c_2$ such that
\beg{align}
c_1\|\pp_t g\|_{L^p_q(T)}\leq \|(\pp_t+b_1\cdot\nn )Jg \|_{L^p_q(T)}\leq c_2\|\pp_t g\|_{L^p_q(T)},~g\in  \cW^{2,p }_{1,q}(T).\label{ppt-bnn-1}
\end{align}
\end{lem}

\beg{proof}
(1) For any $h\in C^2_c(\R^d)$ and  $j\in \{0,1,2\}$, it follows from  \eqref{sup-ps} that
\beg{align}\label{nn-g-ps}
\nn^{j} (J_th)(x) & = h( \ps^{-1}(t,x))\1_{[j=0]}+\sum_{i=1}^d\left[( \nn^j(\ps^{-1})^i)\pp_i h( \ps^{-1})\right](t,x)\1_{[j\geq 1]}\nonumber\\
&\qquad +\left[(\nn \ps^{-1})^*\nn^2 h(  \ps^{-1})(\nn \ps^{-1})\right](t,x)\1_{[j=2]}. 
\end{align}
Then 
\beg{align*}
\int_{\R^d}|\nn^{j} (J_th)(x)|^p\d x &= \int_{\R^d}|\nn^{j} (J_t h)( \ps(t,x))|^p|\det \nn \ps(t,x)|\d x \\
&\leq C_{T,\|\nn b_1\|_\infty,p,j} \sum_{k=0}^j\int_{\R^d}|\nn^k h( x)|^p\d x.
\end{align*}
Since $C^2_c(\R^d)$ is dense in $W^{j,p}(\R^d)$, it is easy to see that $J_t$ is a bounded operator on $W^{j,p}(\R^d)$ with $\sup_{t\in [0,T]}\|J_t\|_{W^{j,p}}<\infty$ for any $j\in\{0,1,2\}$ .  Moreover, \eqref{nn-g-ps} holds for any $h\in W^{j,p}(\R^d)$. Since for any $g\in W^{2,p}_q(T)$, $(Jg)(t,x)=(J_tg(t,\cdot))(x)$. Then   $J$ is a bounded operator in $W^{2,p}_q(T)$. It is easy to see that $J_t$ and $J$ are invertible, and
$$(J^{-1}g)(t,x)=g(t,\ps(t,x)),~g \in W^{2,p}_q(T).$$ 
It is easy to see that $J^{-1}$ is a bounded linear operator.  By \eqref{nn-g-ps}, \eqref{inv-inv} and that for any $g\in W^{2,p}_q(T)$, $g(t,\cdot)\in W^{2,p}(\R^d)$ for a.e. $t\in [0,T]$,  \eqref{nn-u-v-1} follows from the approximation argument. 

(2) For any $g\in C^1_b([0,T],W^{2,p}(\R^d))$, it follows from Lemma \ref{lem-ps} that
\beg{align}\label{add-pptJ-com}
(\pp_t Jg)(t,\cdot)&=\pp_t(g(t,\ps^{-1}(t,\cdot)))\nonumber \\
&=\pp_t g(t,\ps^{-1}(t,\cdot))+\left\<(\nn g)(t,\ps^{-1}(t,\cdot)),\ff {\d \ps^{-1}} {\d t}(t,\cdot)\right\>\nonumber \\
&=\pp_t g(t,\ps^{-1}(t,\cdot))-\left\<(\nn g)(t,\ps^{-1}(t,\cdot)),(\nn\ps)^{-1}(t,\ps^{-1}(t,\cdot))b_1(t,\cdot)\right\>\nonumber \\
&=\pp_t g(t,\ps^{-1}(t,\cdot))-\left\<(\nn g)(t,\ps^{-1}(t,\cdot)),(\nn\ps^{-1})(t,\cdot)b_1(t,\cdot)\right\>\nonumber \\
&=\pp_t g(t,\ps^{-1}(t,\cdot))-\left\<(\nn\ps^{-1})^*(t,\cdot)(\nn g)(t,\ps^{-1}(t,\cdot)),b_1(t,\cdot)\right\>\nonumber \\
&=(J\pp_tg)(t,\cdot)-\left\<(\nn Jg(t,\cdot)),  b_1(t,\cdot)\right\>.
\end{align}
Combining this with \eqref{lin-gro},  we have that 
\beg{align*}
\| Jg\|_{W^{p,w}_{1,q}}^q&\leq 2^{q-1}\left(\int_0^T\|J\pp_t g(t,\cdot)\|_{L^p_w}^q\d t+\int_0^T \|(\nn Jg)\cdot b_1(t,\cdot) \|_{L^p_w}^q\d t\right)\\
& \leq 2^{q-1}\left\{\int_0^T\|J\pp_t g(t,\cdot)\|_{L^p }^q\d t+\int_0^T \|\nn Jg (t,\cdot) \|_{L^p }^q\d t\right\}\\
& \leq 2^{q-1}C_T\left(\int_0^T\| \pp_t g(t,\cdot)\|_{L^p }^q\d t+\int_0^T \|  g (t,\cdot) \|_{W^{1,p}}^q\d t\right).
\end{align*}
Hence, 
\beg{align*}
\| Jg\|_{W^{p,w}_{1,q}(T)}\leq C_T\left(\| g\|_{W^{p,1}_{1,q}(T)}+\|g\|_{W^{2,p}_{q}(T)}\right)=C_T \|g\|_{\cW^{2,p}_{1,q}(T)}. 
\end{align*}
Combining this with that $C^1_b([0,T],W^{2,p}(\R^d))$ is dense in $\cW^{2,p}_{1,q}(T)$ and $J$ is bounded in $W^{2,p}_q(T)$, we have that  $J$ is bounded from $\cW^{2,p}_{1,q}(T)$ to $\cW^{2,p,w}_{1,q}(T)$. Consequently, \eqref{pu-pv-bnu} follows from \eqref{add-pptJ-com}.  Moreover,  
\beg{align}\label{ppt-bnnu}
\int_{\R^d} |(\pp_t +b_1\cdot \nn)Jg|^p(t,x)\d x&= \int_{\R^d} |(\pp_t+b_1\cdot \nn) Jg|^p(t,\ps(t,x))|\det(\nn \ps(t,x))|\d x\nonumber\\
&= \int_{\R^d} |\pp_t g|^p(t,x)|\det(\nn \ps(t,x))|\d x,
\end{align}
and 
\beg{align*}
\int_{\R^d} |\pp_t g|^p(t,x) \d x&=\int_{\R^d} |\pp_t g |^p(t,\ps^{-1}(t,x))|\det(\nn\ps^{-1})(t,x)|\d x\\
&=\int_{\R^d} |(\pp_t  -b_1\cdot \nn)Jg|^p(t,x)|\det(\nn\ps^{-1})(t,x)|\d x.
\end{align*}
Hence, \eqref{ppt-bnn-1} follows from \eqref{sup-ps}.

Finally, for any $g\in \cW^{2,p}_{1,q}(T)$, \eqref{ppt-bnn-1} implies that $Jg\in \cW^{2,p,w }_{1,q,b_1}(T)$. Due to the open mapping theorem, to prove that $J$ is a homeomorphism, it suffices to show that $J$ is a surjective from $\cW^{2,p}_{1,q}(T)$ to $\cW^{2,p,w}_{1,q,b_1}(T)$.  

For any $g\in W_{1,q}([0,T],L_w^p(\R^d))$, we have by \cite[Theorem 1.16]{Bar} that 
\beg{align}\label{deri-1}
&\varlimsup_{\ep\ra 0}\left\|\ff {g(t+\ep,\ps (t+\ep,\cdot))-g(t,\ps (t+\ep,\cdot))} {\ep}-\pp_t g(t,\ps (t+\ep,\cdot))\right\|_{L_w^p}\nonumber\\
&\qquad \leq C_T\varlimsup_{\ep\ra 0}\left\|\ff {g(t+\ep,\cdot)-g(t,\cdot)} {\ep}-\pp_t g(t,\cdot)\right\|_{L_w^p}\nonumber\\
&\qquad =0,~\mbox{a.e.}~t\in [0,T].
\end{align}
Since $\pp_t g(t,\cdot)\in L_w^p(\R^d)$, for any $\de>0$, there exists $g_\de \in C^2_c(\R^d)$ such that $\|g_\de (\cdot)-\pp_t g(t,\cdot)\|_{L_w^p}<\de$. Then
\beg{align}\label{ppt-app}
&\|\pp_t g(t,\ps (t+\ep,\cdot))-\pp_t g(t,\ps (t,\cdot))\|_{L_w^p}\nonumber\\
&\qquad \leq 2C_T\de+\|g_\de (\ps (t+\ep,\cdot))-g_\de(\ps (t ,\cdot))\|_{L_w^p}.
\end{align}
Since that $\nn g_\de$ has compact support and that  there exists $C_T>0$ such that
$$\sup_{t\in [0,T]}|\ps(t,x)|\leq C_{T,K_0}(1+|x|),$$
there exists $M>0$  such that such that
$$\bigcup_{\th\leq 1,|\ep|\leq 1}\supp\left( \nn g_\de(\ps (t ,\cdot)+\th(\ps (t+\ep,\cdot)-\ps(t ,\cdot)))\right)\subset B(0,M),$$
where $B(0,M)$ is the ball center at $0$ with radius $M$. Then, by the mean value theorem and Lemma \ref{lem-ps},  we have that
\beg{align*}
&\left|g_\de (\ps (t+\ep,x))-g_\de(\ps (t ,x))\right|\\
&\qquad =\left|\int_0^1\< \nn g_\de( \th\ps (t+\ep,x)+(1+\th)\ps (t ,x)),\ps (t+\ep,x)-\ps (t ,x)\>\d \th\right|\\
&\qquad \leq C_{\de,T}\left(\int_0^1 \left|\ps (t+\ep,x)-\ps (t ,x)\right|\d \th\right)\1_{B(0,M)}(x)\\
&\qquad =  C_{\de,T}\left( \left|\int_t^{t+\ep} b_1(r,\ps(r,x))\d r\right|\right)\1_{B(0,M)}(x)\\
&\qquad \leq C_{\de,T }\left( \int_t^{t+\ep}C_T(1+|x|)\d r  \right)\1_{B(0,M)}(x)\\
&\qquad \leq \ep C_{\de,T }(1+|x|)\1_{B(0,M)}(x),
\end{align*}
where we have used \eqref{inv} and the linear growth of $b_1$ in the last second inequality. Then we obtain that
\beg{align*}
\|g_\de (\ps^{-1}(t+\ep,\cdot))-g_\de(\ps^{-1}(t ,\cdot))\|_{L_w^p}\leq  C_{\de,T}\ep.
\end{align*}
Putting this into \eqref{ppt-app}, we  prove  indeed that
\beg{align*}
\lim_{\ep\ra 0}\|\pp_t g(t,\ps (t+\ep,\cdot))-\pp_t g(t,\ps (t,\cdot))\|_{L_w^p}=0.
\end{align*}
Combining this with \eqref{deri-1} and that $g(t,\cdot)\in W^{2,p}(\R^d)$ a.e. $t\in [0,T]$,  we have  in $L^p_w(\R^d)$  
\beg{align*}
&\lim_{\ep\ra 0}\ff 1 {\ep}\left( g(t+\ep,\ps (t+\ep,\cdot))-g(t,\ps (t,\cdot))\right)\\
&\qquad =\lim_{\ep\ra 0}\ff 1 {\ep}\left( g(t+\ep,\ps (t+\ep,\cdot))- g(t,\ps^{-1}(t+\ep,\cdot))\right)\\
&\qquad \qquad +\lim_{\ep\ra 0}\ff 1 {\ep}\left( g(t,\ps (t+\ep,\cdot))-g(t,\ps (t,\cdot))\right)\\
&\qquad =\pp_t g(t,\ps (t,\cdot))+\<\nn g(t,\ps (t,\cdot)),  b_1(t,\ps(t,\cdot))\>.
\end{align*}
Setting $J^{-1}g(t,x)=g(t,\ps(t,x))$, we obtain  
$$\pp_t J^{-1}g(t,\cdot)=\left(J^{-1}\left(\pp_t+b_1\cdot\nn\right) g\right),~g\in \cW^{2,p,w }_{1,q,b_1}(T).$$
Hence, $J^{-1}g\in W_{1,q}([0,T],L^p(\R^d))$. Combining this with that $J$ is a homeomorphism on $W^{2,p}_q(T)$,  we have that  $J$ is a surjective to $\cW^{2,p,w}_{1,q,b_1}(T)$.

\end{proof}

\bigskip

\noindent{\bf{\emph{Proof of Theorem \ref{lem-lplq}:~}}}

Let $v$ be the solution of \eqref{equ-pa-v}. By a direct calculus, $u=Jv$ satisfies \eqref{equ-pa-1} in the weak sense. Then Combining Remark \ref{rem-0}, Lemma \ref{lem-v} and Lemma \ref{lem-J}, we can prove  all the assertions directly by the interpolation of the Sobolev spaces except the estimation for $\|u\|_{C^{\al}_{q_2}(t,T)}$ in \eqref{add-la-be}.  However, by  \eqref{nn-u-v-1}, we have that $u(t,x)=v(t,\ps^{-1}(t,x))$, and 
$$\nn u(t,x)=[\nn\ps^{-1}(t,x)]^* \nn v(t,\ps^{-1}(t,x)),$$
if $v(t,\cdot)\in C^1(\R^d)$. Combining this with Lemma \ref{lem-ps}, there is $C_T>0$ such that  
$$\|u(s,\cdot)\|_{C^{\al}}\leq C_T\|v(s,\cdot)\|_{C^{\al}},~s\in [0,T],$$
which yields that $\|u\|_{C^{\al}_{q_2}(t,T)}\leq C_T\|v\|_{C^{\al}_{q_2}(t,T)}$. Hence, \eqref{add-la-be} holds for $u$.



\section{Krylov's estimate}
Let $X_t$ satisfy the following equation
\beg{align}\label{equ-X-0}
X_t=X_0+\int_0^t (b_1+b_2+b_0)(s,X_s)\d s+\int_0^t\si(s,X_s)\d W_s+\int_0^t\xi(s)\d s,
\end{align}
where $\xi(t)$ is an $\sF_t$-adapted process. We  investigate Krylov's estimate for $X_t$ in this section.  Let $a=\ff 1 2\si\si^*$. By using \eqref{equ-pa-1} and \eqref{equ-pa-v} with $c\equiv 0$, we  can prove the following Krylov's estimate for $X_t$. 

\beg{thm}\label{kry-00}
Let $p,q,p_1,q_1\in (1,\infty)$ with   
$$\ff d p +\ff 2 q<2,\qquad\ff d {p_1}+\ff 2 {q_1}<1.$$
Assume {\bf (H1)}-{\bf (H3)} and $b_0 \in L^{p_1}_{q_1}(T)$. Let $\ta$ be a stopping time and $0\leq t_0<t_1\leq T$.   Then for any $f\in L^p_q(T)$, there exist  a positive constant $C $ depending on $p,d,q,T,K_0,\ka_1,\ka_2$, $\|b_0\|_{L^{p_1}_{q_1}(T)}$, $\|\nn b_1\|_{T,\infty}$,$\|b_2\|_{T,\infty}$ such that   
\beg{align*}
&\E\left(\int_{t_0\we\ta}^{t_1\we\ta}f(s,X_s)\d s\Big|\sF_{t_0 }\right)\leq C\left\{1+\left(\E\left(\int_{t_0\we\ta}^{t_1\we\ta}|\xi(s)|^2\d s\Big| \sF_{t_0 }\right)\right)^{\ff 1 2}\right\}^{ k_{p,q} }\|f\|_{L^p_q(t_0,t_1) },
\end{align*}
where $k_{p,q}$ is the smallest integer that greater that $\log_2\left(\ff 2 { 2-\ff d p-\ff 2 q }\right)$.
\end{thm}

\beg{rem}
If $\xi=0$, under the assumption of Theorem \ref{kry-00} and $\nn\si\in L^p_q(T)$,  it follows from \cite[Theorem 1.1]{ZhangX11} that  \eqref{equ-X-0} has a unique local strong solution. It can be proved that the solution of \eqref{equ-X-0} is nonexplosive, see \eqref{add-UU}, \eqref{add-U-U-1} and \eqref{add-U-1} in the proof of this theorem. We use this Krylov's estimate to establish Harnack inequalities in Section 4.
\end{rem}

\subsection{Proofs of Theorem \ref{kry-00}}

To prove Theorem \ref{kry-00}, we need the following two lemmas. 

\beg{lem}\label{kry-1}
Assume {\bf (H1)}, {\bf (H2')}, {\bf (H3)} and $p,q\in (1,\infty)$. Let $\ta$ be a stopping time and $0\leq t_0\leq t_1\leq T$.\\
(1) If  $b_0\equiv 0$, then for any $f\in L^p_q(T)$ with   $\ff d p+\ff 2 q<2$, there exist positive constants  $C$ depending on $p,d,q,T$, $\|\nn b_1\|_\infty$, $\|b_2\|_\infty$, $\ka_1,\ka_2$  such that
\beg{align*}
\E\left(\int_{t_0\we\ta}^{t_1\we\ta}f(s, X_s)\d s\Big|\sF_{t_0 }\right) \leq   C\left\{1+ \left( \E\left(\int_{t_0\we\ta}^{t_1\we\ta}|\xi(s)|^2\d s\Big| \sF_{t_0 }\right)^{\ff 1 2}\right) \right\}^{k_{p,q}}\|f\|_{L^p_q(t_0,t_1) }.
\end{align*}
Moreover, $\ff d p+\ff 2 q<1$, $|\xi(s)|^2$ can be replaced by $|\xi(s)|$ in the above inequality.\\
(2)  If   $\ff d p+\ff 2 q<1$  and $b_0 \in L^p_q(T)$, then for any $f\in L^p_q(T)$, there exists a positive constant $C $ depending on $p,d,q,T,\|\nn b_1\|_\infty,\|b_2\|_\infty,\ka_1,\ka_2$ such that   
\beg{align*}
&\E\left(\int_{t_0\we\ta}^{t_1\we\ta}f(s,X_s)\d s\Big|\sF_{t_0 }\right)\\
&\qquad \leq   C\left(1+\E\left(\int_{t_0\we\ta}^{t_1\we\ta}|\xi(s)|\d s\Big| \sF_{t_0 }\right)\right)\e^{C\|b_0\|^q_{L^p_q(t_0,t_1)}}\|f\|_{L^p_q(t_0,t_1) }.
\end{align*}
\end{lem}

By using Krylov's estimate in Lemma \ref{kry-1}, we can establish a  generalized It\^o's  formula for $u\in \cW^{2,p,w }_{1,q,\tld b_1}(T)$, where $\tld b_1$ satisfying {\bf (H2')}  can be different from $b_1$. Let  $\tld J$ be defined as $J$ with $b_1$ replaced by $\tld b_1$. Then we have the following  lemma.

\beg{lem}\label{ITO}
Let $ X_t$ be a solution of \eqref{equ-X-0} and $  u\in  \cW^{2,p }_{1,q,\tld b_1}(T)$ with $\ff d p+\ff 2 q<1$. Then $\P$-a.s. 
\beg{align*}
&u(t,X_t)=u(t_0,X_{t_0})+\int_{t_0}^t\pp_su(s,X_s)\d s+\int_0^t\<(b_1+b_2+b_0)(s,X_s)+\xi(s),\nn u(s,X_s)\>\d s\\
&+\ff 1 2 \int_0^t \trac(\si^*\si(s,X_s)\nn^2 u(s,X_s))\d s+\int_0^t\<\nn u(s,X_s),\si(s,X_s)\d W_s\>,~0\leq t_0\leq t_1\leq T.
\end{align*}
\end{lem}

The proofs of Lemma \ref{kry-1} and Lemma \ref{ITO} will be given in the next subsection. Now, we  prove Theorem \ref{kry-00} by using these lemmas.

Since $b_0\in L^{p_1}_{q_1}$ with $\ff d {p_1} +\ff 2 {q_1}<1$ and Theorem \ref{lem-lplq}, we can set that $u=(u^1,\cdots,u^d)$ is a solution of following parabolic system on $[0,t_1]$ with $a=\ff 1 2\si\si^*$:
\beg{align*}
\pp_t u^i+L_tu^i = -b_0^i+\la u^i,~u^i(t_1,\cdot)=0,~i=1,\cdots,d.
\end{align*}
Then $u\in \cW^{2,p_1,w}_{1,q_1, b_1}(t_1)$ due to Theorem \ref{lem-lplq} again. Let $U(t,x)=x+u(t,x)$. By Lemma \ref{ITO}, we have
\beg{align}\label{add-UU}
\d U(t,X_t) &= (b_1(t,X_t)+b_2(t,X_t)+b_0(t,X_t)+\xi(t))\d t+\left(\pp_t u(t,X_t)+L_t u(t,X_t)\right)\d t\nonumber\\
&\qquad +\left(I+\nn u(t,X_t)\right)\si(t,X_t)\d W_t+ \nn u(t,X_t) \xi(t) \d t\nonumber\\
&=(b_1+b_2+\la u)(t,X_t)\d t+\left(I+\nn u(t,X_t)\right)\si(t,X_t)\d W_t\nonumber\\
&\qquad +\left(I+\nn u(t,X_t)\right) \xi(t)\d t.
\end{align}
By  \eqref{add-la-be}, for large enough $\la$, we have that $\|\nn u\|_{t_1,\infty} < 1$. Then $U(s,\cdot):\R^d\ra\R^d$ is a diffeomorphism on $\R^d$ for any  $0\leq s\leq t_1$, and 
\beg{align}\label{add-U-U-1}
 \|\nn U^{-1} \|_{t_1,\infty}+\|\nn U \|_{t_1,\infty} <\infty.
\end{align}
Thus, letting $Y_t=U(t,X_t)$, we have
\beg{align}\label{add-U-1}
\d Y_t & =\left(b_1(t,U^{-1}(t,Y_t)+b_2(t, U^{-1}(t,Y_t))+\la u(t,U^{-1}(t,Y_t))\right)\d t\nonumber\\
&\quad +\left(I+\nn u(t,X_t)\right) \xi(t)\d t+\left(I+\nn u(t,U^{-1}(t,Y_t))\right)\si(t,U^{-1}(t,Y_t))\d W_t.
\end{align}
Since $ b_1(t,U^{-1}(t,x))$ satisfies (H2) and $(b_2+\la u)(t, U^{-1}(t,x))$ satisfies (H3), it follows from Remark \ref{rem-0} that we can apply (1) of Lemma \ref{kry-1} to $Y_t$. Taking into account that
\beg{align*}
\int_{\R^d}\left|f(s,U^{-1}(s,x))\right|^p\d x&=\int_{\R^d}\left|f(s,y)\right|^p|\det(\nn U^{-1}(s,y))|\d x\\
&\leq C\int_{\R^d}\left|f(s,y)\right|^p\d y,~f\in L^p_q(t_0,t_1)\\
|\left(I+\nn u(t,X_t)\right)\xi(t)|&\leq 2|\xi(t)|,~t\in (t_0,t_1).
\end{align*}
we have
\beg{align*}
&\E\left(\int_{t_0\we\ta}^{t_1\we\ta}f(s,X_s)\d s\Big|\sF_{t_0 }\right) = \E\left(\int_{t_0\we\ta}^{t_1\we\ta}f(s,U^{-1}(s,Y_s))\d s\Big|\sF_{t_0 }\right) \\
&\qquad \leq   C\left\{1+  \left(\E\left(\int_{t_0\we\ta}^{t_1\we\ta}| \xi(s)|^2\d s\Big| \sF_{t_0 }\right)\right)^{\ff 1 2} \right\}^{k_{p,q}}\|f(\cdot, U^{-1}(\cdot,\cdot))\|_{L^p_q(t_0,t_1)}\\
&\qquad \leq C\left\{1+  \left(\E\left(\int_{t_0\we\ta}^{t_1\we\ta}|\xi(s)|^2\d s\Big| \sF_{t_0 }\right) \right)^{\ff 1 2}\right\}^{k_{p,q}} \|f\|_{L^p_q(t_0,t_1) }.
\end{align*}
Therefore,  we complete the proof.

\subsection{Proofs of Lemma \ref{kry-1} and Lemma \ref{ITO}}
\noindent{\bf{\emph{Proof of Lemma \ref{kry-1}:~}}}

(1) Let $r=d+1$ and $f\in L^p_q(T)\cap L^r_r(T)$ with any $p,q\in(1,+\infty)$. By Theorem \ref{lem-lplq},  the following equation has a unique solution
\beg{align*}
(\pp_t+L_t)u=f ,~u(t_1,\cdot)=0.
\end{align*}
Let $\ps$ be the solution of \eqref{char-lin} with $T$ replaced by $t_1$ and $J$ be the mapping induced by $\ps^{-1}$ as in Lemma \ref{lem-J}. Then $v:=J^{-1}u$ satisfies \eqref{equ-pa-v} with $\bar c=0$ and $T$ replaced by $t_1$. Let $\rh$ be a non-negative smooth function on $\R^{d+1}$ with compact support in the unit ball centre at zero and $\int_{\R^{d+1}}\rh(t,x)\d t\d x=1$. Set $\rh_n(t,x)=n^{d+1}\rh(nt,nx)$. Extending $v$ by zero for $t\geq t_1$ and by $v(0,\cdot)$ for $t<0$. We define 
\beg{align*}
v_n(t,x) &=\int_{\R^{d+1}}v(t-s,x-y)\rh_n(s,y)\d s\d y,\\
\bar f_n&=(\pp_t+\bar L_t )v_n.
\end{align*}
Since $b_0=0$ and $\bar b_2$ is bounded, according to the proof of  \cite[Theorem 2.1]{ZhangX11}, there is 
$$\lim_{n\ra+\infty}\|\bar f_n-\bar f\|_{L^r_r(t_1)}=0.$$
Let $f_n=J\bar f_n$. Then, owing to Lemma \ref{lem-J}, we have 
$$\lim_{n\ra+\infty}\| f_n- f\|_{L^r_r(t_1)}=\lim_{n\ra+\infty}\| J\bar f_n- J\bar f\|_{L^r_r(t_1)}=0.$$
This, combining with \cite[Lemma 3.1]{GM}, yields that 
\beg{align}\label{add-fn-f}
\lim_{n\ra\infty}\E\left(\int_0^{\ta\we t_1}\left|f_n(s,X_s)-f(s,X_s)\right|\d s\right)\leq C\lim_{n\ra\infty}\|f_n-f\|_{L^r_r(t_1)}=0.
\end{align}
By \eqref{inv} and {\bf (H2')},  $\ps^{-1}(t,x)$ is absolutely continuous in $t$ and $\ff {\d \ps^{-1}} {\d t}(t,x)$ is continuous in $x$ uniformly w.r.t $t\in [0,t_1]$.  Then $u_n:=Jv_n$ is absolutely continuous in $t$ and $(\pp_t u_n)(t,x)$ is continuous in $x$ uniformly w.r.t $t\in [0,t_1]$. Moreover, due to Lemma \ref{lem-ps} and that $$\|\nn v_n\|_{t_1,\infty}+\|\nn^2 v_n\|_{t_1,\infty}<\infty,$$
we have that $u_n\in C_b^{0,2}([0,t_1]\times\R^d)$.  Then, we can apply the It\^o formula to $u_n(t,X_t)$ and obtain that 
\beg{align}\label{ineq-un-Ito}
u_n(t,X_t)&=u_n(0,X_0)+\int_0^t\left(\pp_s+L_s \right)u_n (s,X_s)\d s+\int_0^t\xi(s)\cdot\nn u_n(s,X_s)\d s\nonumber\\
&\qquad +\int_0^t\<\nn u_n(s,X_s),\si(s,X_s)\d W_s\>\nonumber\\
&=u_n(0,X_0)+\int_0^tf_n(s,X_s)\d s+\int_0^t\xi(s)\cdot\nn u_n(s,X_s)\d s \nonumber\\
&\qquad +\int_0^t\<\nn u_n(s,X_s),\si(s,X_s)\d W_s\>.
\end{align}
Since $\|\nn u_n\|_{t_1,\infty}<\infty$, Doob's optional theorem yields 
\beg{align*}
\E\left(\int_{t_0\we\ta}^{t_1\we\ta}\<\nn u_n(s,X_s),\si(s,X_s)\d W_s\>\Big|\sF_{t_0}\right)=0.
\end{align*}
Then
\beg{align}\label{ineq-un-Ito-2}
&\E\left(u_n(t_1\we\ta,X_{t_1\we\ta})-u_n(t_0\we\ta,X_{t_0\we\ta})\Big| \sF_{t_0 }\right)\nonumber\\
&\qquad \geq \E\left(\int_{t_0\we\ta}^{t_1\we \ta} f_n(s,X_s)\d s\Big| \sF_{t_0 }\right)\nonumber\\
&\qquad\qquad-\sup_{s\in [t_0,t_1]}\|\nn u_n(s)\|_\infty\E\left(\int_{t_0\we\ta}^{t_1\we\ta}|\xi(s)|\d s\Big|\sF_{t_0 }\right).
\end{align}

We first assume that $\ff d p+\ff 2 q<1$. Then, by Lemma \ref{lem-v} (or \eqref{add-la-be} in Theorem \ref{lem-lplq})  and Lemma \ref{lem-J}, we have that
\beg{align*}
\sup_{n\geq 1,t\in[t_0,t_1]}\left(\|u_n(t,\cdot)\|_{\infty}+\|\nn u_n \|_{\infty}\right)&=\sup_{n\geq 1,t\in[t_0,t_1]}\left(\|(\nn\ps^{-1})^*\nn v_n(t,\cdot) \|_{ \infty}+\|v_n(t,\cdot) \|_{ \infty}\right)\\
&\leq C_1\sup_{n\geq 1,t\in[t_0,t_1]} \left(\|\nn v_n(t,\cdot)\|_{ \infty}+\|v_n (t,\cdot)\|_{ \infty}\right)\\
&\leq  C_1\sup_{ t\in[t_0,t_1]}\left(\|\nn v(t,\cdot) \|_{ \infty}+\|v(t,\cdot)  \|_{\infty}\right)\\
& \leq  C_2   \|f \|_{L^p_q(t_0, t_1 )},
\end{align*}
where $C_1>0$ is a constant depending on $t_1,\|\nn b_1\|_{t_1,\infty}$ and $C_2>0$ is a constant depending on $t_1,p,q$, $\|\nn b_1\|_{t_1,\infty}$, $\|b_2\|_{t_1,\infty},\ka_1,\ka_2$. 
Putting this into \eqref{ineq-un-Ito-2} and taking $n\ra\infty$, we have
\beg{align*}
&\E\left(\int_{t_0\we\ta}^{t_1\we\ta}f_n(s,X_s)\d s\Big|\sF_{t_0 }\right)\nonumber\\
&\qquad \leq C\left(1+\E\left(\int_{t_0\we\ta}^{t_1\we\ta}|\xi(s)|\d s\Big| \sF_{t_0 }\right)\right)\|f_n\|_{L^p_q(t_0,t_1)},
\end{align*}
which together with \eqref{add-fn-f} implies by  that
\beg{align}\label{Krylov-app-n}
&\E\left(\int_{t_0\we\ta}^{t_1\we\ta}f(s,X_s)\d s\Big|\sF_{t_0 }\right)=\lim_{n\ra\infty}\E\left(\int_{t_0\we\ta}^{t_1\we\ta}f_n(s,X_s)\d s\Big|\sF_{t_0 }\right)\nonumber\\
&\qquad  \leq   2\sup_{n\geq 1,t\in [t_0,t_1]}\|u_n(t,\cdot)\|_\infty+  \E\left(\int_{t_0\we\ta}^{t_1\we\ta}|\xi(s)|\d s\Big| \sF_{t_0 }\right) \|f\|_{L^p_q(t_0,t_1)}\nonumber\\
&\qquad  \leq   C\left(1+  \E\left(\int_{t_0\we\ta}^{t_1\we\ta}|\xi(s)|\d s\Big| \sF_{t_0 }\right)\right)\|f\|_{L^p_q(t_0,t_1)}
\end{align}
where $C$ depends on $d,p,q$, $t_1,t_0$, $\|\nn b_1\|_{t_1,\infty},\|b_2\|_{t_1,\infty}$, $\ka_1,\ka_2$. 

Next, we introduce an iteration. Let 
$$\et=\left(\E\left(\int_{t_0\we\ta}^{t_1\we\ta}|\xi(s)|^2\d s\Big| \sF_{t_0 }\right)\right)^{\ff 1 2},\qquad \ga_k=2(1-2^{-k}),~k\in\N.$$
We assume that for any $k\geq 1$ and any $f\in L^{p }_{q }(t_0,t_1)$ with $p ,q \in (1,+\infty)$ and $\ff d {p } +\ff 2 {q }<\ga_k$, there is  
\beg{align}\label{krylov-k}
&\E\left(\int_{t_0\we\ta}^{t_1\we\ta}f(s,X_s)\d s\Big|\sF_{t_0 }\right)\leq C\left(1+ \et\right)^k\|f\|_{L^p_q(t_0,t_1)}.
\end{align} 
Then we prove that \eqref{krylov-k} holds for any $f\in L^{p }_{q }(t_0,t_1)$  with $\ff d {p } +\ff 2 {q }<\ga_{k+1}$ and $k$ replaced by $k+1$ in the right hand side.  Since $\lim_{k\ra+\infty} \ga_k=2$ and that \eqref{krylov-k} holds for $k=1$, the assertion of (1) follows from the induction. 

Let $f_n,\bar f, u_n, u,v_n,v,\ps$ be defined as above, and let $p'=\ff {\ga_{k+1}} {\ga_k}p$, $q'=\ff {\ga_{k+1}} {\ga_k}q$. By \eqref{krylov-k}, we have that 
\beg{align*}
\E\left(\int_{t_0\we \ta}^{t_1\we\ta}\xi(s)\cdot\nn u_n(s,X_s)\d s\Big|\sF_{t_0 } \right)&\leq \et\left(\E\left(\int_{t_0\we \ta}^{t_1\we\ta} |\nn u_n(s,X_s)|^2\d s\Big|\sF_{t_0 } \right)\right)^{\ff 1 2}\\
&\leq C \et (1+\et )^k  \| |\nn u_n|^2 \|_{L^{p'}_{q'}(t_0,t_1)}^{\ff 1 2}\\
&=C\et(1+\et)^{k} \|  \nn u_n \|_{L^{2p'}_{2q'}(t_0,t_1)}.
\end{align*}
Since 
\beg{align*}
1+\ff d {2p'}+\ff 2 {2q'}-\ff d p-\ff 2 q& =1+ \ff {\ga_{k}} {2\ga_{k+1}}\left(\ff d p+\ff 2 q\right)-\left(\ff d p+\ff 2 q\right)\\
&= 1-\ff {2\ga_{k+1}-\ga_k} {2\ga_{k+1}}\left(\ff d p+\ff 2 q\right)\\
&> 1-\ga_{k+1}+\ff 1 2\ga_k\\
&=0,
\end{align*}
it follows from \eqref{add-la-be}, Lemma \ref{lem-ps} and $v_n=v*\rh_n$ that 
\beg{align*}
\varlimsup_{n\ra+\infty}\| \nn u_n \|_{L^{2p'}_{2q'}(t_0,t_1)}&\leq C_1\varlimsup_{n\ra+\infty}\|\nn v_n\|_{L^{2p'}_{2q'}(t_0,t_1)}= C_1\|\nn v\|_{L^{2p'}_{2q'}(t_0,t_1)}\\
&\leq C_2 \|\bar f\|_{L^{p}_q(t_0,t_1)}\leq \tld C_2 \| f\|_{L^{p}_q(t_0,t_1)}.
\end{align*}
Putting this into \eqref{ineq-un-Ito}, as \eqref{ineq-un-Ito-2} and \eqref{Krylov-app-n}, we have that 
\beg{align*} 
&\E\left(\int_{t_0\we\ta}^{t_1\we\ta}f(s,X_s)\d s\Big|\sF_{t_0 }\right)=\lim_{n\ra\infty}\E\left(\int_{t_0\we\ta}^{t_1\we\ta}f_n(s,X_s)\d s\Big|\sF_{t_0 }\right) \\
&\qquad  \leq   2\sup_{n\geq 1,t\in [t_0,t_1]}\|u_n(t,\cdot)\|_\infty+ C(1+ \et)^k\et \|f\|_{L^p_q(t_0,t_1)} \\
&\qquad  \leq   C\left(1+   \et\right)^{k+1} \|f\|_{L^p_q(t_0,t_1)}.
\end{align*}
By the definition of $\ga_k$, it is clear that $k>\log_2(2/(2-\ff d p-\ff 2 q))$. Then we obtain $k_{p,q}$. 

(2) We investigate  the case that $b_0\in L^p_q(T)$ with $\ff d p+\ff 2 q<1$. Set $m>0$, $M>0$, and let
$$\ta_M=\inf\{t\geq 0~|~\int_0^t|b_0(s,X_s)|\d s\geq M\}.$$
It is clear that 
\beg{align*}
X_t& =X_0+\int_0^t \left(b_1(s,X_s)\d s+ b_2(s,X_s)+b_0(s,X_s)\1_{[|b_0(s,X_s)|<m]}\right)\d s\\
&\qquad +\int_0^t \left(b_0(s,X_s)\1_{[|b_0(s,X_s)|\geq m]}+\xi(s)\right) \d s+\int_0^t\si(s,X_s)\d W_s,~t>0.
\end{align*}
Let $L_t^{[m]}=L_t-(b_0\1_{[|b_0|\geq m]})\cdot\nn$. Consider 
$$(\pp_t+L_t^{[m]})u=f,~u(t_1,\cdot)=0.$$
Let $J$ be defined as in the proof of (1) and $v=J^{-1} u$. Arguing as in \cite[Theorem 2.2]{ZhangX11} and using the mapping $J$ and $v$, we can prove that 
\beg{align}\label{add-kry00}
&\E\left(\int_{t_0\we\ta\we\ta_M}^{t_1\we\ta\we\ta_M}f(s,X_s)\d s\Big|\sF_{t_0 }\right)\nonumber\\
&\qquad  \leq   C\left(1+\E\left(\int_{t_0\we\ta\we\ta_M}^{t_1\we\ta\we\ta_M}\left(\xi(s)+|b_0(s,X_s)|\1_{[|b_0(s,X_s)|\geq m]}\right)\d s\Big| \sF_{t_0 }\right)\right)\nonumber\\
&\qquad \qquad \times\e^{C\|b_0\1_{[|b_0|\leq m]}\|_{L^p_q(t_0,t_1)}^q}\|f\|_{L^p_q(t_0,t_1) },
\end{align}
where $C$ is a positive constant as claimed  in (1) of Theorem \ref{lem-lplq}  which is independent of $m$. For fixed $M>0$,  we have that
\beg{align*}
\lim_{m\ra\infty}\E \int_{t_0\we\ta\we\ta_M}^{t_1\we\ta\we\ta_M}|b_0(s,X_s)|\1_{[|b_0(s,X_s)|\geq m]}\d s=0.
\end{align*}
Thus, taking  $m\ra+\infty$ first and then $M\ra+\infty$ in \eqref{add-kry00}, we complete the proof of the second assertion of this lemma.

\bigskip

\noindent{\bf{\emph{Proof of Lemma \ref{ITO}:~}}}

Let  $\rh$ be a non-negative smooth function on $\R^{d }$ with compact support in the unit ball centre at zero and $\int_{\R^{d }}\rh( x) \d x=1$, and let $\rh_n(x)=n^{d}\rh(nx)$. For any $n\geq 1$, define 
$$u_n(t,x)=\int_{\R^d}u(t,y)\rh_n( x-y)\d y,~x\in\R^d.$$

We now divide the proof into two steps.

{\emph{Step (i): we are going to show that the conclusion holds for $u_n(t,x)$.}}

For any $m\in\N$,  Let  $t_k=\ff {k t} {m}$, $k=0,1,\cdots,m$. Then
\beg{align}\label{pre-Ito}
u_n(t,X_t)-u_n(0,X_0) & =\sum_{k=0}^{m-1}\left(u_n(t_{k+1},X_{t_{k+1}})-u_n(t_{k},X_{t_{k+1}})\right)\nonumber\\
&\qquad+\sum_{k=0}^{m-1}\left(u_n(t_{k },X_{t_{k+1}})-u_n(t_{k},X_{t_{k}})\right)\nonumber\\
& =: I_{1,m}+I_{2,m}. 
\end{align}

We first study  $I_{1,m}$. Let
$$\pp_t u_n(t,x)=\int_{\R^d}\pp_t u(t, y)\rh_n(x-y)\d y.$$
Then for any $M>0$, we have by Jessen's inequality  that
\beg{align}\label{sup-ppun}
\sup_{|x|\leq M}\left|\pp_t u_n(t,x)\right| & \leq \sup_{|x|\leq M}\left(\int_{\R^d} |\pp_t u(t, y)|^p\rh_n(x-y)\d y\right)^{\ff 1 p}\nonumber\\
& \leq \left(\sup_{|x|\leq M}(1+|y|^2)^{\ff 1 2}\rh_n^{\ff 1 p}(x-y)\right) \left\|\pp_t u(t,\cdot)\right\|_{L^p_w}\nonumber\\
& \equiv C_{M,n,p}\left\|\pp_t u(t,\cdot)\right\|_{L^p_w}.
\end{align}
This implies that $\sup_{|x|\leq M}\left|\pp_t u_n(\cdot,x)\right|\in L^q([0,T])$ since $u\in \cW^{2,p,w }_{1,q,\tld b_1}(T)$. 
Moreover,
\beg{align*}
&\sup_{|x|\leq M}\left|\ff {u_n(t+\ep,x)-u_n(t,x)} {\ep}-\pp_t u_n(t,x)\right|\\&\qquad=\sup_{|x|\leq M}\left|\int_{\R^d}\left(\ff {u(t+\ep,y)-u(t,y)} {\ep}-\pp_t u(t,y)\right)\rh_n(x-y)\d y\right|\\
&\qquad \leq \sup_{|x|\leq M}\left(\int_{\R^d}\left|\ff {u(t+\ep,y)-u(t,y)} {\ep}-\pp_t u(t,y)\right|^p\rh_n(x-y)\d y\right)^{\ff 1 p}\\
&\qquad \leq \left(\sup_{|x|\leq M,y\in\R^d}(1+|y|^2)^{\ff 1 2}\rh_n^{\ff 1 p}(x-y)\right)\left\|\ff {u(t+\ep,\cdot)-u(t,\cdot)} {\ep}-\pp_t u(t,\cdot)\right\|_{L^p_w}\\
&\qquad \equiv C_{M,p,n} \left\|\ff {u(t+\ep,\cdot)-u(t,\cdot)} {\ep}-\pp_t u(t,\cdot)\right\|_{L^p_w}.
\end{align*}
Because  $u\in \cW^{2,p,w }_{1,q,\tld b_1}(T)$, $u$ is absolutely continuous in $L^p_w(\R^d)$.  Then for a.e. $t\in [0,T]$, we have
\beg{align*}
&\varlimsup_{\ep\ra 0}\sup_{|x|\leq M}\left|\ff {u_n(t+\ep,x)-u_n(t,x)} {\ep}-\pp_t u_n(t,x)\right|\\
&\qquad \leq  C_{M,p,n} \varlimsup_{\ep\ra 0}\left\|\ff {u(t+\ep,\cdot)-u(t,\cdot)} {\ep}-\pp_t u(t,\cdot)\right\|_{L^p_w}\\
&\qquad =0.
\end{align*}
Hence,
\beg{align*}
I_{1,m}& =  \sum_{k=0}^{m-1} \int_0^t(\pp_s u_n)(s,X_{t_{k+1}})\1_{[t_k\leq s<t_{k+1}]}(s)\d s\\
& =  \int_0^t\sum_{k=0}^{m-1}(\pp_s u_n)(s,X_{t_{k+1}})\1_{[t_k\leq s<t_{k+1}]}(s)\d s.
\end{align*}
By the H\"older inequality, for any $|x_1|\leq M$ and $|x_2|\leq M$, we have
\beg{align*}
&\left|\pp_t u_n(t,x_1)-\pp_t u_n(t,x_2)\right|\\
&\qquad \leq \int_{\R^d}|\pp_t u(t,y)|\int_0^1 |\nn \rh_n(x_2+\th (x_1-x_2)-y)|\d \th\d y|x_1-x_2|\\
&\qquad \leq   |x_1-x_2|\|\pp_t u(t,\cdot)\|_{L^p_w}\\
&\qquad\qquad \times\left(\int_{\R^d}\int_0^1 |\nn \rh_n(x_2+\th (x_1-x_2)-y)|^{\ff p {p-1}}(1+|y|^2)^{\ff p {2(p-1)}}\d \th\d y\right)^{\ff {p-1} p}\\
&\qquad  \leq C_{M,n,p}|x_1-x_2|\|\pp_t u(t,\cdot)\|_{L^p_w},~\mbox{a.e.}~t\in [0,T].
\end{align*} 
Combining this  with the continuity of the path of $X_t$, we have that for a.e. $s\in [0,T]$ 
$$\lim_{m\ra\infty}\sum_{k=0}^{m-1}(\pp_t u_n)(s,X_{t_{k+1}})\1_{[t_k\leq s<t_{k+1}]}(s)= \pp_s u_n(s,X_s).$$
Taking into account that  $\P$-a.s.  $\tld M:=\sup_{s\in [0,T]}|X_s|<\infty$ and the inequality \eqref{sup-ppun}, we have that
\beg{align*}
\left|\sum_{k=0}^{m-1}(\pp_t u_n)(s,X_{t_{k+1}})\1_{[t_k\leq s<t_{k+1}]}(s)\right|&\leq   \sum_{k=0}^{m-1} |\sup_{0\leq k\leq m-1}(\pp_t u_n)(s,X_{t_{k+1}})|\1_{[t_k\leq s<t_{k+1}]}(s)\\
& \leq C_{\tld M,n,p}\sum_{k=0}^{m-1}\left\|\pp_t u(s,\cdot)\right\|_{L^p_w}\1_{[t_k\leq s<t_{k+1}]}(s)\\
&= C_{\tld M,n,p} \left\|\pp_t u(s,\cdot)\right\|_{L^p_w}. 
\end{align*}
Hence, it follows from the dominated convergence theorem that  $\P$-a.s. 
\beg{align*}
\lim_{m\ra\infty}I_{1,m}=\int_0^t \pp_t u_n(s,X_s)\d s,~t\in [0,T].
\end{align*}

Next, we investigate $I_{2,m}$. Since $u\in \cW^{2,p,w }_{1,q,\tld b_1}(T)$, we have $\tld J^{-1}u\in \cW^{2,p}_{1,q}(T)$. Then it follows from \cite[Lemma 10.2]{KR} and \eqref{inv} that $u\in C(\R\times\R^d)$. For any $t\in [0,T]$, by the definition of $u_n$, \cite[Lemma 10.2]{KR} and Lemma \ref{lem-J}, we have $u_n(t,\cdot)\in C^2(\R^d)$ and for any $j\in\{0,1,2\}$
\beg{align}\label{sup-nnun}
\|\nn^j u_n\|_{T,\infty} & \leq \|u\|_{T,\infty} \int_{\R^d}|\nn^j\rh_n(x-y)|\d y\nonumber\\
&= \|\tld J^{-1} u\|_{T,\infty}  \int_{\R^d}|\nn^j\rh_n( y)|\d y\nonumber\\
&\leq C_{n,j}C_{T,p,q,d}\left(\|\pp_t \tld J^{-1} u\|_{L^p_q(T)}+T\|\tld J^{-1} u\|_{W^{2,p}_q(T)}\right).
\end{align}
For $j\in \{0,1,2\}$, any $t,s\in [0,T]$ and   $M>0$,
\beg{align}\label{nnu-nnu}
&\sup_{|x|\leq M}\left\|\nn^{j}u_n(t,x)-\nn^{j}u_n(s,x)\right\| \nonumber\\
&\qquad\leq  \sup_{|x|\leq M}\int_{\R^d}|u(t,y)-u(s,y)||\nn^{j}\rh_n(x-y)|\d y\nonumber\\
&\qquad \leq |t-s|\sup_{|x|\leq M}\int_0^1 \int_{\R^d}|\pp_t u(s+\th(t-s),y)||\nn^j\rh_n(x-y)|\d y \d \th\nonumber\\
&\qquad\leq C_{M,n,q,j}|t-s|\left(\int_0^1\|\pp_t u(s+\th(t-s),\cdot)\|^q_{L^p_w}\d \th\right)^{\ff 1 q}\nonumber\\
&\qquad\leq C_{M,n,q,j}\|u\|_{W^{p,w}_{1,q}(T)}|t-s|^{\ff {q-1} q}.
\end{align}

By It\^o's formula, we have
\beg{align*}
I_{2,m}& =\sum_{k=0}^{m-1}\left(\int_{t_k}^{t_{k+1}}\<(b_1+b_2+b_0)(s,X_s)+\xi(s),\nn u_n(t_k,X_s)\>\d s\right.\\
& \qquad +\ff 1 2\int_{t_k}^{t_{k+1}}\trac\left(\si\si^*(s,X_s)\nn^2 u_n(t_k,X_s)\right)\d s\\
&\qquad +\left. \int_{t_k}^{t_{k+1}}\<\nn u_n(t_k,X_s),\si(s,X_s)\d W_s\>\right)\\
& =\int_{0}^{t}\<(b_1+b_2+b_0)(s,X_s)+\xi(s),\sum_{k=0}^{m-1}\nn u_n(t_k,X_s)\1_{[t_k\leq s<t_{k+1}]}(s)\>\d s \\
& \qquad +\ff 1 2\int_{0}^{t}\sum_{k=0}^{m-1}\trac\left(\si\si^*(s,X_s)\nn^2 u_n(t_k,X_s)\right)\1_{[t_k\leq s<t_{k+1}]}(s)\d s\\
&\qquad + \int_{0}^{t}\<\sum_{k=0}^{m-1}\nn u_n(t_k,X_s)\1_{[t_k\leq s<t_{k+1}]}(s),\si(s,X_s)\d W_s\> .
\end{align*}
It follows from \eqref{nnu-nnu} that
\beg{align*} 
&\left|\sum_{k=0}^{m-1}\left\{\left(\nn u_n(t_k,X_s)-\nn u_n(s,X_s)\right)+\left(\nn^2 u_n(t_k,X_s)-\nn^2 u_n(s,X_s)\right)\right\}\1_{[t_k\leq s<t_{k+1}]}(s)\right|\\
&\qquad \leq C_{M,n,q}\left(\ff t m\right)^{\ff {q-1} q},
\end{align*}
where $M=\sup_{s\in[0,t]} |X_s|$. Then 
\beg{align*}
&\lim_{m\ra\infty}\sum_{k=0}^{m-1}\left(\nn u_n(t_k,X_s)+\nn^2 u_n(t_k,X_s)\right)\1_{[t_k\leq s<t_{k+1}]}(s)\\
&\qquad  =\nn u_n(s,X_s)+\nn^2 u_n(s,X_s).
\end{align*}
Since $\int_0^T\left( |b_1(s,X_s)|+|b_2(s,X_s)|+|b_0(s,X_s)|+|\xi(s)|\right)\d s<\infty$ and $\|\si \|_\infty<\infty$, it follows from \eqref{sup-nnun} and the dominated convergence theorem that $\P$-a.s.
\beg{align*}
&\lim_{m\ra\infty}\left(\int_{0}^{t}\<(b_1+b_2+b_0)(s,X_s)+\xi(s),\sum_{k=0}^{m-1}\nn u_n(t_k,X_s)\1_{[t_k\leq s<t_{k+1}]}(s)\>\d s \right.\\
& \qquad \qquad \left. +\ff 1 2\int_{0}^{t}\sum_{k=0}^{m-1}\trac\left(\si\si^*(s,X_s)\nn^2 u_n(t_k,X_s)\right)\1_{[t_k\leq s<t_{k+1}]}(s)\d s\right)\\
&\qquad = \left(\int_{0}^{t}\<(b_1+b_2+b_0)(s,X_s)+\xi(s), \nn u_n(s,X_s) \>\d s \right.\\
& \qquad \qquad \left. +\ff 1 2\int_{0}^{t} \trac\left(\si\si^*(s,X_s)\nn^2 u_n(s,X_s)\right) \d s\right),~t\in [0,T].
\end{align*}
Let $\tld \ta_M=\inf\{s>0~|~|X_s|\geq M\}$. Then \eqref{nnu-nnu} yields that 
\beg{align*}
&\varlimsup_{m\ra\infty}\E\int_0^{t\we\tld\ta_M}\left|\si^*(s,X_s)\sum_{k=0}^{m-1}\left(\nn u_n(t_k,X_s)-\nn u_n(s,X_s)\right)\1_{[t_k\leq s<t_{k+1}]}(s) \right|^2\d s\\
&\qquad \leq \lim_{m\ra\infty}C_{M,n}\|\si\|^2\left(\ff t m\right)^{\ff {2(q-1)} q}t\\
&\qquad =0.
\end{align*}
Hence, $\P$-a.s. 
\beg{align*}
\lim_{m\ra\infty}I_{2,m}&=\int_{0}^{t}\<(b_1+b_2+b_0)(s,X_s)+\xi(s), \nn u_n(s,X_s) \>\d s \\
& \qquad +\ff 1 2\int_{0}^{t} \trac\left(\si\si^*(s,X_s)\nn^2 u_n(s,X_s)\right) \d s\\
&\qquad + \int_{0}^{t}\< \nn u_n(s,X_s) ,\si(s,X_s)\d W_s\>. 
\end{align*}

By letting $m\ra\infty$ in \eqref{pre-Ito}, we obtain that  It\^o's formula for $u_n(t,x)$:
\beg{align}\label{ito-un}
u_n(t,X_t)-u_n(0,X_0) &=\int_0^t \pp_t u_n(s,X_s)\d s\nonumber\\
&\qquad +\int_{0}^{t}\<(b_1+b_2+b_0)(s,X_s)+\xi(s), \nn u_n(s,X_s) \>\d s \nonumber\\
& \qquad +\ff 1 2\int_{0}^{t} \trac\left(\si\si^*(s,X_s)\nn^2 u_n(s,X_s)\right) \d s\nonumber\\
&\qquad + \int_{0}^{t}\< \nn u_n(s,X_s) ,\si(s,X_s)\d W_s\>\nonumber\\
&=: I_1(n)+I_2(n)+I_3(n)+I_4(n).
\end{align}

{\emph{Step (ii): we shall complete the proof of this lemma by using the approximation argument.}}

It is clear that for any $j\in\{0,1\}$, we have $\|\nn^j u_n\|_{T,\infty}\leq \|\nn^j u\|_{T,\infty}$. It follows from \eqref{nn-u-v-1}, $u\in \cW^{2,p,w }_{1,q,\tld b_1}(T)$, \cite[Lemma 10.2]{KR} and Lemma \ref{lem-J} that 
\beg{align*}
\sup_{t\in [0,T]}\left|u(t,x)-u(t,y)\right|&\leq C_{T,\|\nn\tld b_1\|_{T,\infty},\| \tld J^{-1} u\|_{\cW^{p,1}_{1,q}(T)}}|x-y|,\\
\sup_{t\in [0,T]}\left|\nn u(t,x)-\nn u(t,y)\right|&\leq C_{T,\|\nn\tld b_1\|_{T,\infty}}\sup_{t\in [0,T]}\left|(\tld J\nn (\tld J^{-1}u))(t,x)-(\tld J\nn (\tld J^{-1}u))(t,y)\right|\\
&\qquad +C_{T,\|\nn\tld b_1\|_{T,\infty}}\|\tld J\nn\tld J^{-1}u\|_{T,\infty} |x-y|\\
&\leq C_{T,\|\nn\tld b_1\|_{T,\infty},\| \tld J^{-1} u\|_{\cW^{p,1}_{1,q}(T)}}\left|x-y\right|^{\de'}\\
&\qquad +C_{T,\|\nn\tld b_1\|_{T,\infty},\| \tld J^{-1} u\|_{\cW^{p,1}_{1,q}(T)}} |x-y|,
\end{align*} 
with some $\de'<1-\ff d p-\ff 2 q$. Then for any $j\in\{0,1\}$
\beg{align*}
&\| \nn^j u_n - \nn^j u\|_{T,\infty}\\
&\qquad=\sup_{x\in\R^d,t\in [0,T]}\int_{\R^d}\left|\nn^j u(t,x-y)-\nn^j u(t,x)\right|\rh_n(y)\d y\\
&\qquad =\sup_{x\in\R^d,t\in [0,T]} \int_{\R^d}\left|\nn^j u(t,x-\ff y n)-\nn^j u(t,x)\right|\rh(y)\d y\\
&\qquad \leq \ff {C_{T,\|\nn\tld b_1\|_{T,\infty},\| \tld J^{-1} u\|_{\cW^{p,1}_{1,q}(T)}}\int_{\R^d}\left(|u|\vee |u|^{\de'}\right)\rh(u)\d u} {n^{\de'}},
\end{align*}
which implies $\lim_{m\ra\infty}\|\nn^j u_n-\nn^j u\|_{T,\infty}=0.$
Applying the dominated convergence theorem to  $I_2(n)$ and $I_4(n)$, we get $\P$-a.s. (by a subsequence if necessary)
\beg{align*}
\lim_{n\ra\infty}\left(I_2(n)+I_4(n)\right) &=\left(\int_{0}^{t}\<(b_1+b_2+b_0)(s,X_s)+\xi(s), \nn u (s,X_s) \>\right.\\
& \qquad\left.+\int_{0}^{t}\< \nn u (s,X_s) ,\si(s,X_s)\d W_s\>\right.,~t\in [0,T].
\end{align*}
Let $R>0$ and $\ta_R=\inf\{t\in [0,T]~|~\int_0^t|\xi(s)|\d s>R\}$. Then it follows from (2) of Lemma \ref{kry-1} that
\beg{align*}
&\E\int_0^{T\we\ta_R}\|\si\si^*(\nn^2 u_n-\nn^2 u)(s,X_s)\|_{HS}\d s\\
&\qquad\leq \|\si\|^2_\infty\E\int_0^{T\we\ta_R}\|(\nn^2 u_n-\nn^2 u)(s,X_s)\|_{HS}\d s\\
&\qquad \leq C\left(1+\E\int_0^{T\we\ta_R}|\xi(s)|\d s\right)\|\nn^2 u_n-\nn^2 u\|_{L^p_q(T)}.
\end{align*}
It is clear that for any $t\in [0,T]$, 
\beg{align*}
\|\nn^2 u_n(t,\cdot)\|_{L^p }&\leq \|\nn^2 u(t,\cdot)\|_{L^p},\\
\lim_{n\ra\infty}\|\nn^2 u_n(t,\cdot)&-\nn^2 u(t,\cdot)\|_{L^p}=0.
\end{align*}
Then by the dominated convergence theorem that
\beg{align*}
\lim_{n\ra\infty}\int_{0}^T\|\nn^2u_n(t,\cdot)-\nn^2 u(t,\cdot)\|_{L^p}^q\d t=0,
\end{align*}
which implies 
$$\lim_{n\ra\infty}\E\int_0^{T\we\ta_R}\|\si\si^*(\nn^2 u_n-\nn^2 u)(s,X_s)\|_{HS}\d s=0.$$
Thus  there exists a subsequence $u_{n_k}$ such that on $\{\ta_R>T\}$
\beg{align}\label{lim-nn2un}
&\lim_{k\ra\infty}\int_{0}^{t} \trac\left(\si\si^*(s,X_s)\nn^2 u_{n_k}(s,X_s)\right) \d s\nonumber\\
&\qquad=\int_{0}^{t} \trac\left(\si\si^*(s,X_s)\nn^2 u(s,X_s)\right) \d s,~t\in [0,T].
\end{align} 
Since $\lim_{R\ra\infty}\P(\ta_R\leq T)=0$, by Cantor's diagonal argument, there exists a subsequence, which we also denote by $u_{n_k}$, such that \eqref{lim-nn2un} holds $\P$-a.s.

For $I_1(n)$, it follows from the property of convolution that
\beg{align*}
\left|\ff {\pp_t u_n(s,x)-\pp_t u(s,x)} {\sq{1+|x|^2}}\right| & \leq \int_{\R^d}\ff {|\pp_t u(s,x-y)|} {\sq{1+|x-y|^2}}\ff {\left|\sq{1+|x-y|^2}-\sq{1+|x|^2}\right|} {\sq{1+|x|^2}}\rh_n(y)\d y\\
&\qquad +\int_{\R^d}\left|\ff {\pp_t u(s,x-y)} {\sq{1+|x-y|^2}}-\ff {\pp_t u(s,x)} {\sq{1+|x|^2}}\right|\rh_n(y)\d y\\
&\leq \int_{\R^d}\ff {|\pp_t u(s,x-y)|} {\sq{1+|x-y|^2}}\ff {|y|} {\sq{1+|x|^2}}\rh_n(y)\d y\\
&\qquad +\int_{\R^d}\left|\ff {\pp_t u(s,x-y)} {\sq{1+|x-y|^2}}-\ff {\pp_t u(s,x)} {\sq{1+|x|^2}}\right|\rh_n(y)\d y\\
&\leq \ff 1 n\int_{\R^d}\ff {|\pp_t u(s,x-y)|} {\sq{1+|x-y|^2}}\rh_n(y)\d y\\
&\qquad +\int_{\R^d}\left|\ff {\pp_t u(s,x-y)} {\sq{1+|x-y|^2}}-\ff {\pp_t u(s,x)} {\sq{1+|x|^2}}\right|\rh_n(y)\d y\\
&\equiv J_{1,n}(s)+J_{2,n}(s).
\end{align*}
By convolution inequality, 
\beg{align*}
\|J_{1,n}(s)\|_{L^p}&\leq \ff 1 n\|\pp_t u(s,\cdot)\|_{L^p_w},\qquad \|J_{2,n}(s)\|_{L^p}\leq 2\|\pp_t u(s,\cdot)\|_{L^p_w},\\
\varlimsup_{n\ra\infty} \|J_{2,n}(s)\|_{L^p}^p&\leq \varlimsup_{n\ra\infty}\int_{\R^d}\left(\int_{\R^d}\left|\ff {\pp_t u(s,x-\ff y n)} {\sq{1+|x-\ff y n|^2}}-\ff {\pp_t u(s,x)} {\sq{1+|x|^2}}\right|^p\d x\right)\rh(y)\d y\\
&=0.
\end{align*}
Hence  
\beg{align}
\left\|\ff {\pp_t u_n(s,\cdot)-\pp_t u(s,\cdot)} {\sq{1+|\cdot|^2}}\right\|_{L^p}\leq 3\|\pp_t u(s,\cdot)\|_{L^p_w},\label{supn-pps-un}\\
\lim_{n\ra\infty}\left\|\ff {\pp_t u_n(s,\cdot)-\pp_t u(s,\cdot)} {\sq{1+|\cdot|^2}}\right\|_{L^p}=0.\label{pps-un-u}
\end{align}
Let $\tld \ta_R=\inf\{t\in [0,T]~|~|X_t|>R\}$. Then by (2) of Lemma \ref{kry-1}
\beg{align*}
&\E\int_0^{T\we\ta_R\we\tld\ta_R}\left|\pp_t u_n(s,X_s)-\pp_t u(s,X_s)\right|\d s\\
&\qquad \leq \sq{1+R^2}\E\int_0^{T\we\ta_R\we\tld\ta_R}\ff {\left|\pp_t u_n(s,X_s)-\pp_t u(s,X_s)\right|} {\sq{1+|X_s|^2}}\d s\\
&\qquad \leq C\sq{1+R^2}\left(1+\E\int_0^{T\we\ta_R\we\tld\ta_R}|\xi(s)|\d s\right)\|\pp_t u_n-\pp_t u\|_{L^{p,w}_q(T)}.
\end{align*}
It follows from \eqref{supn-pps-un}, \eqref{pps-un-u} and the dominated convergence theorem that 
$$\lim_{n\ra\infty}\|\pp_t u_n-\pp_t u\|_{L^{p,w}_q(T)}=0,$$
which implies that 
$$\lim_{n\ra\infty}\E\int_0^{T\we\ta_R\we\tld\ta_R}\left|\pp_su_n(s,X_s)-\pp_s u(s,X_s)\right|\d s=0.$$
Hence, by the Cantor's diagonal argument, there exists a subsequence, denoted also by $u_n$, such that $\P$-a.s.
$$\lim_{n\ra\infty}\int_0^t\left|\pp_t u_n(s,X_s)-\pp_t u(s,X_s)\right|\d s=0,~t\in [0,T].$$

Combining all these together, we complete the proof by taking $n\ra\infty$ in \eqref{ito-un}.

\section{Applications}

Consider  \eqref{equ-main-0} with $a=\ff 1 2\si\si^*$ satisfying {\bf (H1)} , and  
\beg{description}[align=left, noitemsep]
\item [(H4)] for every $T>0$, $\sup_{t\in [0,T] }|b(t,0)|<\infty$ and  $\|\nn b\|_{T,\infty}<\infty$,   and  there exist $p,q\in (1,\infty)$ with $\ff d p+\ff 2 q<1$ such that $b_0,\nn \si \in L^p_q(T)$. 
\end{description}
By Remark \ref{rem-0}, {\bf (H1)}, {\bf (H4)} and Theorem \ref{kry-00},  we can follow the proofs of \cite[Theorem 1.1]{ZhangX11} or \cite[Theorem 2.1]{XZ} to prove \eqref{equ-main-0} has a unique strong solution.  Let $P_t$ be the associated semigroup generated by $X_t$.  In this section, we shall investigate Harnack inequalities for \eqref{equ-main-0}.

We first establish  log-Harnack inequality following the methodology of \cite{LLW}.
\beg{thm}\label{log-Har}
Assume  {\bf (H1)} and {\bf (H4)}. Then there exists $K_1>0$ such that
\beg{align}\label{log-H}
P_T\log f(x)\leq \log P_T f(y)+\ff {K_1  |x-y|^2} {\ka_1T},~x,y\in\R^d,f\in \sB_b^+(\R^d).
\end{align}
\end{thm}

Fix some $T>0$. By Remark \ref{rem-0}, {\bf (H1)} and {\bf (H4)}, it follows from Theorem  \ref{lem-lplq} that the following system has a  unique solution
\beg{align}\label{equ-pa-v2}
\pp_t \ph^i+\trac(a\nn^2 \ph^i)+(b+b_0)\cdot\nn \ph^i=-b_0^i+\la \ph^i,~\ph^i(T,x)=0,~i=1,\cdots, d,
\end{align}
with $a=\ff 1 2\si\si^*$. Let $\ph =(\ph^1,\cdots,\ph^d)$ and $\Ph_s(x)=x+\ph(s,x)$, $s\in [0,t]$. Since $\ff d p+\ff 2 q<1$, it follows from \eqref{add-la-be}  that we can choose large enough $\la$ such that $\|\nn \ph\|_{T,\infty} <\ff 1 2$. Then 
\beg{align}\label{Ph}
\ff 1 2 |x-y|&\leq |\Ph_s(x)-\Ph_s(y)|\leq \ff 3 2 |x-y|, \\
\ff 2 3 |x-y|&\leq |\Ph_s^{-1}(x)-\Ph_s^{-1}(y)|\leq  2 |x-y|,~s\in [0,T],~x,y\in\R^d.\label{Ph-inv}
\end{align}

By Lemma  \ref{ITO}, we have
\beg{align*}
\d \Ph_s(X_s)&=\left(\pp_s\Ph_s(X_s)+L_s\Ph_s(X_s)\right)\d s+\left(I+\nn \ph(s,X_s)\right)\si(s,X_s)\d W_s\\
&=\left(b(s,X_s)+b_0(s,X_s)\right)\d s+\left(\pp_s \ph(s,X_s)+L_s\ph(s,X_s)\right)\d s\\
&\qquad +\left(I+\nn \ph(s,X_s)\right)\si(s,X_s)\d W_s\\
&=\left(b(s,X_s)+\la \ph(s,X_s)\right)\d s +\left(I+\nn \ph(s,X_s)\right)\si(s,X_s)\d W_s.
\end{align*}
Let $Y_s=\Ph_s(X_s)$ and $X_s=\Ph_s^{-1}(Y_s)$. Define 
\beg{align*}
Z(s,y)=(b+\la u)(s,\Ph_s^{-1}(y)),\qquad \Si(s,y)=\left(I+\nn \ph(s,\Ph_s^{-1}(y))\right)\si\left(s,\Ph_s^{-1}(y)\right).
\end{align*}
Then we transform \eqref{equ-main-0}  to 
\beg{align}\label{equ-y}
\d Y_s= Z(s,Y_s)\d s+\Si(s,Y_s)\d W_s,\qquad Y_0=\Ph_0(X_0).
\end{align}
It is clear that 
$$\nn\Si\in L^p_q(t),\qquad \ff 1 4 \ka_1|x|^2 \leq |\Si^*(s,y)x|^2\leq \ff 9 4 \ka_2|x|^2,~s\in [0,T],~x,y\in\R^d,$$ 
and there exists $\bar K_1>0$ such that 
$$|Z(s,y_1)-Z(s,y_2)|\leq \bar K_1|y_1-y_2|,~y_1,y_2\in \R^d,~s\in [0,T].$$

Let $\cT_sf(x)=\E f(Y_s^x)$ with $Y_0^x\equiv x$. Hence, we have the following log-Harnack inequality for \eqref{equ-y}. Since $\|\nn \Si\| \in L^p_q(T)$ and $\|\nn Z\|_{T,\infty}<\infty$, the proof of this lemma follows from that of  \cite[Proposition 2.1]{LLW} completely.
\beg{lem}
There exists $\tld K_0$ such that for any $f\in \sB_b^+(\R^d)$,
\beg{align*}
\cT_T\log f(y)\leq \log \cT_T f(x)+\ff {\tld K_0 |x-y|^2} {\ka_1 T},~x,y\in\R^d. 
\end{align*}
\end{lem}

Since $P_sf(x)=\E f(X_s^x)=\E f(\Ph_s^{-1}(Y_s^{\Ph_0(x)}))=\cT_s \bar f(\Ph_0(x))$ with $\bar f(\cdot)=f(\Ph_s^{-1}(\cdot))$, Theorem \ref{log-Har} follows from this lemma and \eqref{Ph} directly.


\bigskip

Next, we shall establish the Harnack inequality with power for \eqref{equ-main-0}. Before our detailed discussions, we give some remarks on the Harnack inequality with power for SDEs with irregular coefficients.

\beg{rem}\label{rem-Har}
For $b_0\in  L^p_q(t)$ with $p,q\in (1,\infty)$ satisfying $\ff d p+\ff 2 q<1$,  since $\|\nn\Si\|\in L^p_q(t)$ does not yield that $\nn\Si(t,\cdot)$ is bounded. So \eqref{equ-y} does not fulfill  conditions to derived Harnack inequalities with power in \cite{W11}.  \cite{Shao} established Harnack inequalities with an  extra constant for SDEs whose drift merely satisfies the $L^p$-$L^q$ integral condition.  Since there exists an extra constant, the Harnack inequality established in \cite{Shao} can not yield the strong Feller property of the associated semigroup.

Recently,  the author in \cite{Huang} assume that the non-regular drift $b_0$ satisfies $b_0\in L^p_q(T)$ with $\ff d p+\ff 2 q<1$ and 
\beg{align}\label{C-Huang}
\int_{\R^d}|b_0(t,x+y)-b_0(t,x)|^p\d x\leq K^p(t)|y|^p
\end{align}
with $K\in L^q_{loc}([0,\infty))$. Then Harnack inequalities were derived. However,   the conditions used in \cite{Huang} are not real conditions that allow the drift to be singular in space variable. In fact, given $t\in [0,T]$, $b_0(t,\cdot)\in L^p(\R^d)$ and \eqref{C-Huang} imply by the definition of the Besov space $B_{p,\infty}^\ga$, see \cite[Section 2.5.12]{Tr}, that $b_0(t,\cdot)\in \cap_{\ga<1}B_{p,\infty}^\ga$ and there exists $C>0$ such that 
$$\|b_0(t,\cdot)\|_{B_{p,\infty}^\ga}\leq C(K(t)+1),~\mbox{a.e.}~t\in [0,T].$$
Since $\ff d p+\ff 2 q<1$, there exists  $\ga-\ff d p>0$ for $\ga$ being closed to $1$. By the embedding theorem of Besov space, see \cite[Theorem 6.5.1]{BL}, we have 
$$b_0(t,\cdot)\in B^{\ga-\ff d p}_{\infty,\infty}(\R^d)= C^{\ga-\ff d p}(\R^d),~\ff d p<\ga<1.$$
Then $b_0(t,\cdot)$ is bounded and $(\ga-\ff d p)$-H\"older continuous for any $\ga\in (\ff d p,1)$ and  there exists $C>0$ such that
$$ \|b_0 \|_{ C^{\ga-\ff d p}_q(T)}^q \leq C\int_0^T(K(t)+1)^q\d t.$$
\end{rem}

Our main result on the Harnack inequality with power is the following theorem, which can be applied to SDEs with singular drift without extra regularity assumption and the diffusion coefficient is H\"older continuous with order greater than $\ff 1 2$.

\beg{thm}\label{HHar}
Fix $T>0$. Assume {\bf (H1)} and {\bf (H4)}, and that there exist $c_T>0$ and $\be>0$ such that
\beg{align}\label{Ho-si}
\|\si(t,x)-\si(t,y)\|_{HS}\leq c_T |x-y|^{\be},~x,y\in\R^d. 
\end{align} 
(1) If $\be\geq \ff 1 2$,  then there are $K_T,K_{ 1,T}>0$ such that for any $\ga>1+4\left(\sq{1+\ff {\la_1^2} {8\la_2}}-1\right)^{-1}$, we have the following Harnack inequality with extra constant
\beg{align}\label{HH-1}
( P_t f )^\ga(y)\leq  P_t f^\ga(x)\exp\left\{  T K_{1,T}+\ff { K_{ T} |x-y|^{2} } { (1-e^{-K_TT})}\right\},~f\in \sB_b^+(\R^d).  
\end{align}
(2) If  $\ff d p+\ff 2 q<\ff 1 2$ and $\be>\ff 1 2$, then there exists $K_{ T}>0$   such that for any $\ga>(1+\ff {6\sq 2\la_2} {\sq\la_1 \al_0 })^2$ with $\al_0=(1-\ff d p-\ff 2 q)\we\be$, any $\al\in (\ff 1 2, 1-\ff d p-\ff 2 q)\cap (\ff 1 2,\be]$, we have 
\beg{align*}
( P_t f )^\ga(y)\leq  P_t f^\ga(x)\exp\left\{\ff {\sq \ga(\sq \ga-1)K_{ T}(|x-y|^2\vee |x-y|^{2\al})} { 2\de_{\ga,T}( \sq\la_T \al (\sq \ga-1)-2\de_{\ga,T})(1-e^{-K_TT})}\right\}.  
\end{align*}
where $\de_{\ga,T}=\ff {3\la_2} 2\vee \ff {\sq\la_1 \al (\sq \ga-1)} {4\sq 2}$  and $f\in \sB_b^+(\R^d)$.
\end{thm}

The proof of (1) of Theorem \ref{HHar}  just follows  that in \cite{Shao} and Krylov's estimate in Theorem \ref{kry-00} directly,  and we leave it to readers. We focus on establishing Harnack inequality without extra constant. To this aim, we first investigate the following equation
\beg{align}\label{equ-hhX}
\d \hh X_t=\hh b(t,\hh X_t)\d t+\hh\si(t,\hh X_t)\d W_t,
\end{align}
where $\hh b:[0,\infty)\times\R^d\ra\R^d$ and $\hh\si[0,\infty)\times\R^d\ra\R^d\otimes\R^d$  satisfying the following   conditions
\beg{description}[align=left, noitemsep]
\item [(H5)] Fix $T>0$. Assume there exist $\al\in (\ff 1 2 ,1]$ and positive constants $K_T$, $\de_T$, $\la_T$ such that for $x,y\in\R^d,~t\in [0,T]$
\beg{align*}
&2\<\hh b(t,x)-\hh b(t,y),x-y\>+\|\hh\si(t,x)-\hh\si(t,y)\|_{HS}^2\leq K_T|x-y|^2\vee |x-y|^{2\al},\\
 & |(\hh\si(t,x)-\hh\si(t,y))^*(x-y)|\leq \de_T|x-y|,~~~
 \hh\si(t,x)\hh\si^*(t,x) \geq \la_T.
\end{align*}\end{description}
Let $\hh X_t$ solve (\ref{equ-hhX}) with $\hh X_0=x.$  Set $\et(t)=\ff {2\al-\th} {K_T}(1-e^{K_T(t-T)})$ with $\th\in (0,2\al)$ and let $\hh Y_t$
 solve the following equation with $Y_0=y$
 \beg{align}\label{equ-hhY}
\d \hh Y_t=\hh b(t,\hh Y_t)\d t+\hh \si(t,\hh Y_t)\d W_t+\ff {\hh\si(t,\hh Y_t)\hh\si^{-1}(t,\hh X_t)(\hh X_t-\hh Y_t)} {\et(t)(|\hh X_t-\hh Y_t|^{2-2\al}\we 1)}\1_{[0,T)}(t)\d t, 
\end{align}
The following lemma  is crucial to establish the Harnack inequality for \eqref{equ-main-0}.
\beg{lem}\label{Har-pow}
Fix $T>0$. Assume {\bf (H5)}. Suppose that \eqref{equ-hhX} has a unique strong solution and for any $(x,y)\in \R^{2d}$, the martingale solution to the system $(\hh X_t,\hh Y_t)$ with $(\hh X_0,\hh Y_0)=(x,y)$ is well-posed on $[0,T)$. Let  
\beg{align}
R_s&=\exp\left\{-\int_0^s\left\<\ff { \hh\si^{-1}(t,\hh X_t)(\hh X_t-\hh Y_t)} {\et(t)(|\hh X_t-\hh Y_t|^{2-2\al}\we 1)},\d W_t\right\>\right.\nonumber\\
&\qquad\qquad\qquad \left.-\ff 1 2\int_0^s \left|\ff { \hh\si^{-1}(t,\hh X_t)(\hh X_t-\hh Y_t)} {\et(t)(|\hh X_t-\hh Y_t|^{2-2\al}\we 1)}\right|^2\d t\right\},~s\in[0,T].\label{Rs}
\end{align}
Then we have the following two conclusions.
\begin{enumerate}
\item[{\rm (i)}] For $\ga_0=\ff {\la_T\th^2} {8(2\de_T+\sq\la_T\th)\de_T},$ one has 
\beg{align}\label{R1q0}
  \sup_{s\in [0,T]}\E R_s^{1+\ga_0}\leq \exp\left\{\ff {(4\de_T+\sq\la_T\th)\th K(|x-y|^2\vee |x-y|^{2\al})} {16(2\de_T+\sq\la_T\th)(2\al-\th)(1-e^{-K_TT})\de_T^2}\right\}.
\end{align}
 \item[{\rm (ii)}] Let $\hh P_t$ be the  associated  transition semigroup of $\hh X_t$. For any $f\in\sB_b^+(\R^d)$, $x,y\in\R^d$ and $\ga>\left(1+\ff {2\de_T} {\sq\la_T \al }\right)^2$,  we have
\beg{align}
(\hh P_t f )^\ga(y)\leq \hh P_t f^\ga(x)\exp\left\{\ff {\sq \ga(\sq \ga-1)K_T(|x-y|^2\vee |x-y|^{2\al})} {4\de_{\ga,T}( \sq\la_T \al (\sq \ga-1)-4\de_{\ga,T})(1-e^{-K_TT})}\right\},  
\end{align}
where $\de_{\ga,T}=\de_T\vee \ff {\sq\la_T \al (\sq \ga-1)} 4$.
\end{enumerate}
\end{lem}
\beg{proof}
Fix $(x,y)\in \R^{2d}$. Since the martingale solution of the system $(\hh X_t,\hh Y_t)$ is well-posed, there exist a system of process $(\hh X_t,\hh Y_t,W_t)_{t\in [0,T)}$ and a probability space with filtration $(\Om,\sF,\P,\sF_t)_{t\in [0,T)}$ such that $\{W\}_{t\in [0,T)}$ is a Brownian motion w.r.t.  $(\Om,\sF,\P,\sF_t)_{t\in [0,T)}$ and $(\hh X_t,\hh Y_t,W_t)_{t\in [0,T)}$ satisfies \eqref{equ-hhX} and \eqref{equ-hhY}. Let $\ta_n=\inf\{t\in[0,T)~|~|\hh X_t|+|Y_t|\geq n\}$, and let
\beg{align*}
\tld W_t=W_t+\int_0^t \ff {\hh\si^{-1}(s,\hh X_s)(\hh X_s-\hh Y_s)} {\et(s)(|\hh X_s-\hh Y_s|^{ 2-2\al }\we 1)}\d s,~t\in[0,T).
\end{align*}
Then the system \eqref{equ-hhX} and \eqref{equ-hhY} can be rewritten as 
\beg{align*}
\d \hh X_t&=\hh b(t,\hh X_t)\d t+\hh\si(t,\hh X_t)\d \tld W_t-\ff {\hh X_t-\hh Y_t} {\et(t)(|\hh X_t-\hh Y_t|^{2-2\al}\we 1)}\1_{[0,T)}(t)\d t,\\
\d \hh Y_t&=\hh b(t,\hh Y_t)\d t+\hh\si(t,\hh Y_t)\d \tld W_t,
\end{align*}
and it follows from Girsanov's theorem that $\{\tld W_t\}_{0\leq t\leq s\we\ta_n}$ is a Brownian motion under $R_{s\we\ta_n}\P$ for any $s\in [0,T)$.   By  It\^o's formula, 
\beg{align*}
\d |\hh X_t-\hh Y_t|^2&=2\<\hh b(t,\hh X_t)-\hh b(t,\hh Y_t),\hh X_t-\hh Y_t\>\d t-\ff {2|\hh X_t-\hh Y_t|^{2\al}\vee |\hh X_t-\hh Y_t|^2} {\et(t)}\d t\\
&\qquad +2\<\hh X_t-\hh Y_t,(\hh\si(t,\hh X_t)-\hh\si(t,\hh Y_t))\d \tld W_t\>\\
&\qquad +\|\hh\si(t,\hh X_t)-\hh\si(t,\hh Y_t)\|_{HS}^2\d t\\
&\leq K_T|\hh X_t-\hh Y_t|^2\vee |\hh X_t-\hh Y_t|^{2\al}\d t-\ff {2|\hh X_t-\hh Y_t|^{2\al}\vee |\hh X_t-\hh Y_t|^2} {\et(t)}\d t\\
&\qquad +2\<\hh X_t-\hh Y_t,(\hh\si(t,\hh X_t)-\hh\si(t,\hh Y_t))\d \tld W_t\>,\qquad t<s\we\ta_n.
\end{align*}
Then
\beg{align*}
\d \ff {|\hh X_t-\hh Y_t|^2} {\et(t)}&\leq \ff {K_T\et(t)-2} {\et^2(t)}|\hh X_t-\hh Y_t|^{2\al}\vee |\hh X_t-\hh Y_t|^2\d t-\ff {\et'(t)} {\et^2(t)}|\hh X_t-\hh Y_t|^2\d t\\
&\qquad +\ff 2 {\et(t)}\<\hh X_t-\hh Y_t,(\hh\si(t,\hh X_t)-\hh\si(t,\hh Y_t))\d \tld W_t\>\\
& \leq -\ff {\et'(t)+2-K_T\et(t)} {\et^2(t)}|\hh X_t-\hh Y_t|^{2\al}\vee |\hh X_t-\hh Y_t|^2\d t\\
&\qquad +\ff 2 {\et(t)}\<\hh X_t-\hh Y_t,(\hh\si(t,\hh X_t)-\hh\si(t,\hh Y_t))\d \tld W_t\>\\
& \leq -\ff {\th |\hh X_t-\hh Y_t|^{2\al}\vee |\hh X_t-\hh Y_t|^2} {\et^2(t)}\d t\\
&\qquad +\ff 2 {\et(t)}\<\hh X_t-\hh Y_t,(\hh\si(t,\hh X_t)-\hh\si(t,\hh Y_t))\d \tld W_t\>,\qquad t<s\we\ta_n.
\end{align*}
Since $\al>\ff 1 2$, it follows from It\^o's formula that 
\beg{align*}
\d \ff {|\hh X_t-\hh Y_t|^{2\al}} {\et(t)}&\leq -\ff {\et'(t)} {\et^2(t)}|\hh X_t-\hh Y_t|^{2\al}\d t+\ff { \al K_T} { \et(t)}|\hh X_t-\hh Y_t|^{4\al-2}\vee |\hh X_t-\hh Y_t|^{2\al}\d t\\
& \qquad -\ff {2\al} {\et^2(t)}|\hh X_t-\hh Y_t|^{2\al}\vee |\hh X_t-\hh Y_t|^{4\al-2}\d t\\
&\qquad +\left\<\ff {2\al(\hh X_t-\hh Y_t)} {\et(t)|\hh X_t-\hh Y_t|^{2-2\al}},(\hh\si(t,\hh X_t)-\hh\si(t,\hh Y_t))\d \tld W_t\right\>\\
&\leq -\ff {2\al +\et'(t)- \al K_T  \et(t)} {\et^2(t)}|\hh X_t-\hh Y_t|^{2\al}\vee |\hh X_t-\hh Y_t|^{4\al-2}\\
&\qquad +\left\<\ff {2\al (\hh X_t-\hh Y_t)} {\et(t)|\hh X_t-\hh Y_t|^{2-2\al}},(\hh\si(t,\hh X_t)-\hh\si(t,\hh Y_t))\d \tld W_t\right\>\\
&\leq -\ff {\th |\hh X_t-\hh Y_t|^{2\al}\vee |\hh X_t-\hh Y_t|^{4\al-2}} {\et^2(t)}\\
&\qquad +\left\<\ff {2\al (\hh X_t-\hh Y_t)} {\et(t)|\hh X_t-\hh Y_t|^{2-2\al}},(\hh\si(t,\hh X_t)-\hh\si(t,\hh Y_t))\d \tld W_t\right\>,~t<s\we\ta_n.
\end{align*}
Then
\beg{align*}
&\E R_{s\we\ta_n}\left(\int_0^{s\we\ta_n}\ff {|\hh X_t-\hh Y_t|^{2}\vee |\hh X_t-\hh Y_t|^{4\al-2}} {\et^2(t)}\d t\right)\\
&\quad \leq \E R_{s\we\ta_n}\left(\int_0^{s\we\ta_n}\ff {|\hh X_t-\hh Y_t|^{2\al}\vee |\hh X_t-\hh Y_t|^{4\al-2}+|\hh X_t-\hh Y_t|^{2\al}\vee |\hh X_t-\hh Y_t|^{2 }} {\et^2(t)} \d t\right)\\
&\quad \leq \ff {K_T\left(|x-y|^2+|x-y|^{2\al}\right)} {\th (1+\al-\th)(1-e^{-K_TT})},
\end{align*}
which yields that 
\beg{align*}
\E R_{s\we\ta_n}\log R_{s\we\ta_n}\leq \ff {2K_T\left(|x-y|^{2\al}\vee |x-y|^2\right)} {\la_T\th(2\al-\th)(1-e^{-K_TT})},~s\in [0,T),~n\in\N.
\end{align*}
Hence, $\{R_{s\we\ta_n}\}_{s<T,n\in\N}$ is a uniformly integrable martingale.  By martingale convergence theorem and $\ta_n\uparrow T$, $R_{s\we\ta_n}$ can be extended to $T$ such that $\{R_s\}_{s\in [0,T]}$ is a martingale. Moreover, it follows from Fatou's lemma that
\beg{align}\label{RlogR}
\sup_{s\in [0,T]}\E R_s\log R_s\leq \ff {2K_T(|x-y|^{2\al}\vee |x-y|^2)} {\la_T\th(2\al-\th)(1-e^{-K_TT})},
\end{align} 
which also implies that $\{R_s\}_{s\in [0,T]}$ is a uniformly integrable martingale.  Hence, by Girsanov's theorem, $\{\tld W_t\}_{t\in [0,T]}$ is a Brownian motion under $R_T\P$. Moreover, by the pathwise uniqueness of \eqref{equ-hhX}, $\{\hh Y_t\}_{t\in [0,T)}$ can be extended to $T$ such that $\lim_{t\ra T}\hh Y_t=\hh Y_T$, and $\{\hh Y\}_{t\in [0,T]}$ is a weak solution of \eqref{equ-hhX} with starting point $x$ replaced by $y$. Moreover, we have $\hh X_T=\hh Y_T$ $\P$-a.s. since \eqref{RlogR} and $\int_0^T\et^{-2}(t)\d t=\infty$. 

Next, we  prove \eqref{R1q0}. Since
\beg{align*}
&\int_0^s \ff {|\hh X_t-\hh Y_t|^{2}\vee |\hh X_t-\hh Y_t|^{4\al-2}} {\et^2(t)}\d t\leq \ff {2K_T |x-y|^{2}\vee |x-y|^{2\al}} {\th(2\al-\th)(1-e^{-K_TT})}\\
&\qquad +\int_0^s\left\<\ff {2(\hh X_t-\hh Y_t)} {\th\et(t)}+\ff { 2\al (\hh X_t-\hh Y_t)} {\th\et(t)|\hh X_t-\hh Y_t|^{2-2\al}},(\hh\si(t,\hh X_t)-\hh\si(t,\hh Y_t))\d \tld W_t\right\>.
\end{align*} 
Denoting by $\E_{ s,n }$ the expectation w.r.t. $R_{s\we\ta_n}\P$, then for $r>0$
\beg{align*}
&\E_{s,n}\exp\left\{r\int_0^{s\we\ta_n} \ff {|\hh X_t-\hh Y_t|^{2}\vee |\hh X_t-\hh Y_t|^{4\al-2}} {\et^2(t)}\d t\right\}\exp\left\{-\ff {2rK_T(|x-y|^{2\al}\vee |x-y|^2)} {\th(2\al-\th)(1-e^{-K_TT})}\right\}\\
&\leq \E_{s,n}\exp\left\{\ff r {\th}\int_0^{s\we\ta_n}\left\< \left(2  +\ff {2\al} { |\hh X_t-\hh Y_t|^{2-2\al}}\right)\ff {\hh X_t-\hh Y_t} {\et(t)},(\hh\si(t,\hh X_t)-\hh\si(t,\hh Y_t))\d \tld W_t\right\>\right\}\\
&\leq \left(\E_{s,n}\exp\left\{\ff {2r^2} {\th^2}\int_0^{s\we\ta_n}  \left(2  +\ff {2\al} { |\hh X_t-\hh Y_t|^{2-2\al}}\right)^2\ff {|(\hh\si(t,\hh X_t)-\hh\si(t,\hh Y_t))^*(\hh X_t-\hh Y_t)|^2} {\et^2(t)}\d t\right\}\right)^{\ff 1 2}\\
&\leq \left(\E_{s,n}\exp\left\{\ff {2\de_T^2 r^2} {\th^2}\int_0^{s\we\ta_n}  \ff {\left(2|\hh X_t-\hh Y_t|   +2\al |\hh X_t-\hh Y_t|^{2\al-1} \right)^2} {\et^2(t)}\d t\right\}\right)^{\ff 1 2}\\
&\leq \left(\E_{s,n}\exp\left\{\ff {32\de_T^2 r^2} {\th^2}\int_0^{s\we\ta_n}  \ff { |\hh X_t-\hh Y_t|^2\vee |\hh X_t-\hh Y_t|^{4\al-2} } {\et^2(t)}\d t\right\}\right)^{\ff 1 2}.
\end{align*}
By taking $r=\ff {\th^2} {32\de_T^2}$, we have
\beg{align}\label{Rq0}
&\E_{s,n}\exp\left\{\ff {\th^2} {32\de_T^2}\int_0^{s\we\ta_n} \ff {|\hh X_t-\hh Y_t|^{2}\vee |\hh X_t-\hh Y_t|^{4\al-2}} {\et^2(t)}\d t\right\}\nonumber\\
&\qquad \leq \exp\left\{\ff {\th K_T(|x-y|^{2\al}\vee |x-y|^2)} {8\de_T^2(2\al-\th)(1-e^{-K_TT})}\right\}.
\end{align}
By H\"older inequality, we have for any $\ga_1>1$ that
\beg{align*}
\E R_{s\we\ta_n}^{1+\ga_0}&\leq \left(\E_{s,n}\exp\left\{\ff {(\ga_1 \ga_0 +1)\ga_0\ga_1} {2(\ga_1-1)}\int_0^{s\we\ta_n}\ff {|\hh\si^{-1}(t,\hh X_t)(\hh X_t-\hh Y_t)|^2} {\et^2(t)(|\hh X_t-\hh Y_t|^{4-4\al}\we 1)}\d t\right\}\right)^{\ff {\ga_1-1} {\ga_1}}\\
&\leq \left(\E_{s,n}\exp\left\{\ff {(\ga_1 \ga_0 +1)\ga_0\ga_1} {2(\ga_1-1)}\int_0^{s\we\ta_n}\ff {|\hh\si^{-1}(t,\hh X_t)(\hh X_t-\hh Y_t)|^2} {\et^2(t)(|\hh X_t-\hh Y_t|^{4-4\al}\we 1)}\d t\right\}\right)^{\ff {\ga_1-1} {\ga_1}}\\
&\leq \left(\E_{s,n}\exp\left\{\ff {(\ga_1 \ga_0 +1)\ga_0\ga_1} {2(\ga_1-1)\la_T}\int_0^{s\we\ta_n}\ff {| \hh X_t-\hh Y_t |^2\vee | \hh X_t-\hh Y_t |^{4\al-2}} {\et^2(t)}\d t\right\}\right)^{\ff {\ga_1-1} {\ga_1}}.
\end{align*}
Taking $\ga_1=1+\sq{1+\ga_0^{-1}}$ which minimizes $\ff {\ga_1(\ga_1\ga_0+1)} {\ga_1-1}$,  we have 
\beg{align*}
\ff {\ga_1\ga_0(\ga_1\ga_0+1)} {2(\ga_1-1)\sq\la_T}&=\ff {\ga_0(\sq {\ga_0} +\sq{\ga_0+1})^2} {2\la_T}=\ff {\th^2} {32\de_T^2},\\
\ff {\ga_1-1} {\ga_1}&=\ff {4\de_T+\sq\la_T\th} {4\de_T+2\sq\la_T\th}.
\end{align*}
Combining this with \eqref{Rq0}, we have 
\beg{align*}
\E R_{s\we\ta_n}^{1+\ga_0}&\leq \exp\left\{\ff {(\ga_1-1)\th K_T(|x-y|^{2\al}\vee |x-y|^2)} {8\ga_1\de_T^2(2\al-\th)(1-e^{-K_TT})}\right\}\\
&=\exp\left\{\ff {(4\de_T+\sq\la_T\th) \th K_T(|x-y|^{2\al}\vee |x-y|^2)} {16(2\de_T+\sq\la_T\th) \de_T^2(2\al-\th)(1-e^{-K_TT})}\right\}.
\end{align*}
Letting $n\ra\infty$, we get \eqref{R1q0}.

For any $\ga>(1+\ff {2\de_T} {\sq\la_T \al })^2$, we  set $\th=\ff {4\de_T} {\sq\la_T(\sq \ga-1)}$. Then $\th<2\al$ and $\ff 1 {\ga-1}=\ga_0$. Consequently, 
\beg{align*}
&\sup_{s\in [0,T]}\left(\E R_{s }^{\ff \ga {\ga-1}}\right)^{\ga-1}=\sup_{s\in [0,T]}\left(\E R_{s }^{1+\ga_0}\right)^{\ga-1}\\
&\qquad \leq \exp\left\{\ff {(\ga-1)(4\de_T+\sq\la_T\th) \th K_T(|x-y|^{2\al}\vee |x-y|^2)} {16(2\de_T+\sq\la_T\th) \de_T^2(2\al-\th)(1-e^{-K_TT})}\right\}\\
&\qquad = \exp\left\{\ff {\sq \ga(\sq \ga-1) K_T(|x-y|^{2\al}\vee |x-y|^2)} {2\de_T[ 2\al (\sq \ga-1)-4\de_T](1-e^{-K_TT})}\right\}.
\end{align*}
Therefore, 
\beg{align*}
(\hh P_Tf)^\ga(y)&=(\E R_T f(\hh Y_T))^\ga\leq (\E R_T^{\ff \ga {\ga-1}})^{\ga-1}\E f^\ga(\hh Y_T)=(\E R_T^{\ff \ga {\ga-1}})^{\ga-1}\E f^\ga(\hh X_T)\\
&\leq \hh P_Tf^\ga(x)\exp\left\{\ff {\sq \ga(\sq \ga-1) K_T(|x-y|^{2\al}\vee |x-y|^2)} {4\de_T[ \al (\sq \ga-1)-2\de_T](1-e^{-K_TT})}\right\}.
\end{align*}
It is clear that this inequality also holds with $\de_T$ replaced by $\de_{\ga,T}$.

\end{proof}

\noindent{\bf{\emph{Proof of Theorem \ref{HHar}:~}}}

Let $\ph$ be given by \eqref{equ-pa-v2} with $\|\nn\ph\|_{T,\infty}<\ff 1 2$, $\Ph_t(x)=x+\ph(t,x)$ and $Y_t=\Ph_t(X_2)$. By \eqref{add-la-be}, for any $\de\in (0,1-\ff d p-\ff 2 q)$, we have that
$$\sup_{t\in [0,T]}\|\nn\ph(t,x)-\nn\ph(t,y)\|\leq C_T|x-y|^\de.$$
Then 
$$\ff {\ka_1} 2|h|^2\leq  \left|[\si^*(I+\nn\ph)^*](t,\Ph^{-1}_t(x)) h\right|\leq \ff 3 2\ka_2|h|^2,~x,h\in\R^d,t\in[0,T],$$ 
and there exists $C_{T}>0$ such that for any $t\in [0,T]$
$$\|(I+\nn\ph)\si(t,\Ph^{-1}_t(x))-(I+\nn\ph)\si(t,\Ph^{-1}_t(y))\|_{HS}\leq C_T|x-y|^{\de}\vee |x-y|^\be.$$
Setting 
$$\hh b(t,x)=b(t,\Ph_t^{-1}(x))+\la\ph(t,\Ph^{-1}(x)),~\hh\si(t,x)=(I+\nn\ph)\si(t,\Ph^{-1}_t(x)),$$
it is clear that the conditions for $\hh b$ and $\hh\si$ in Lemma \ref{Har-pow} holds with $\de_T=\ff 3 2\la_2$, $\la_T=\ff {\la_1} 2$ and any  $\al\in (\ff 1 2, 1-\ff d p-\ff 2 q)\cap (0,\be]$. Hence, the Harnack inequality with power follows if we prove the well-posedness of the martingale solution to the system $(\hh X_t,\hh Y_t)$ with  $\hh b$ and $\hh\si$ defined as above. 

To get the well-posedness, according to \cite[Corollary 10.1.2]{SV}, we only need to investigate the following system $(\hh X_t^{(n)},\hh Y_t^{(n)})$ with $n\in\N$:
\beg{align*}
\d \hh X_t^{(n)}&=\hh b(t,\hh X_t^{(n)})\d t+\hh\si(t,\hh X_t^{(n)})\d W_t,\\
\d \hh Y_t^{(n)}&=\hh b(t,\hh Y_t^{(n)})\d t+\hh\si(t,\hh Y_t^{(n)})\d   W_t\\
&\qquad +\ff {\hh\si(t,\hh Y_t^{(n)})\hh\si^{-1}(t,\hh X_t^{(n)})} {\et(t)}\pi_n\left(\ff {\hh X_t^{(n)}-\hh Y_t^{(n)}} {|\hh X_t^{(n)}-\hh Y_t^{(n)}|^{2-2\al}\we 1}\right)\1_{[0,T)}(t)\d t,
\end{align*}
where $\pi_n:\R^d\ra \R^d$ defined as follows
$$\pi_n(y)=y\1_{[|y|<n]}+n\ff {y} {|y|}\1_{[|y|\geq n]},~y\in\R^d.$$
Since $\pi_n$ is a bounded function and $\hh\si$ is bounded and non-generated, the well-posedness of $(\hh X_t^{(n)},\hh Y_t^{(n)})$ can be investigated via the following system with any Brownian motion $\{W_t\}_{t\geq 0}$ and Girsanov's theorem:
\beg{align}\label{hhXn}
\d \hh X_t^{(n)}&=\hh b(t,\hh X_t^{(n)})\d t+\hh\si(t,\hh X_t^{(n)})\d W_t\nonumber\\
&\qquad -\et^{-1}(t)\pi_n\left(\ff {\hh X_t^{(n)}-\hh Y_t^{(n)}} {|\hh X_t^{(n)}-\hh Y_t^{(n)}|^{2-2\al}\we 1}\right)\1_{[0,T)}(t)\d t,\\
\d \hh Y_t^{(n)}&=\hh b(t,\hh Y_t^{(n)})\d t+\hh\si(t,\hh Y_t^{(n)})\d   W_t.\label{hhYn} 
\end{align} 
By  {\bf (H1)}, {\bf (H4)}, \eqref{hhYn} has a unique strong solution. Next, we prove \eqref{hhXn} has a pathwise unique solution. For any two solutions of \eqref{hhXn}, say $\hh X_t^{(n),1}$ and $\hh X_t^{(n),2}$ with the same initial value and  $\hh Y_t^{(n)}$, there exists $C>0$ such that
\beg{align*}
\d |\hh X_t^{(n),1}-\hh X_t^{(n),2}|^2&\leq C|\hh X_t^{(n),1}-\hh X_t^{(n),2}|^2\d t+\|\si(t,\hat X_t^{(n),1})-\si(t,\hat X_t^{(n),2})\|_{HS}^2\d t\\
&\qquad +2\< \hh X_t^{(n),1}-\hh X_t^{(n),2} ,(\hh\si(t,\hh X_t^{(n),1})-\hh\si\hh X_t^{(n),2})\d W_t\>,
\end{align*}
where we have used the following inequality
\beg{align*}
\left\<\pi_n\left(\ff y {|y|^{2-2\al}\we1}\right)-\pi_n\left(\ff y {|y|^{2-2\al}\we1}\right),y_1-y_2\right\>\geq 0.
\end{align*}
Since {\bf (H1)}, {\bf (H4)} and $\ph\in W^{2,p}_q(T)$, we have $\hh\si\in W^{1,p}_q(T)$. Hence, by \cite[Lemma 5.4]{ZhangX11},
\beg{align*}
&\|\hh\si(t,\hh X_t^{(n),1})-\hh\si(t,\hh X_t^{(n),2})\|_{HS}\\
&\qquad\leq \left( |\cM(\nn \hh\si(t,\cdot))(\hh X_t^{(n),1})|+|\cM(\nn \hh\si(t,\cdot))(\hh X_t^{(n),2})|\right)|\hh X_t^{(n),1}-\hh X_t^{(n),2}|,
\end{align*}
where $\cM$ is the Hardy-Littlewood maximal function defined by 
$$(\cM f)(x)=\sup_{r\in (0,\infty)}\ff 1 {|B_r|}\int_{B_r}f(x+y)\d y$$
with $B_r=\{x\in\R^d~|~|x|\leq r\}$ and $f$ is a locally bounded function. Consequently, 
\beg{align}\label{path-ito}
\d |\hh X_t^{(n),1}-\hh X_t^{(n),2}|^2&\leq C |\hh X_t^{(n),1}-\hh X_t^{(n),2}|^2\left(1+|\cM(\nn \hh\si(t,\cdot))(\hh X_t^{(n),1})|^2\right.\nonumber\\
&\qquad\qquad\left.+|\cM(\nn \hh\si(t,\cdot))(\hh X_t^{(n),2})|^2\right)\d t\nonumber\\
& \qquad +2\< \hh X_t^{(n),1}-\hh X_t^{(n),2} ,(\hh\si(t,\hh X_t^{(n),1})-\hh\si\hh X_t^{(n),2})\d W_t\>.
\end{align}
Since $\pi_n$ is a bounded function, for any $s<T$,
$$\tld W_r=W_r-\int_0^r\et^{-1}(t)\hh\si^{-1}(t,\hh X_t^{(n),1})\pi_n\left(\ff {\hh X_t^{(n),1}-\hh Y_t^{(n)}} {|\hh X_t^{(n),1}-\hh Y_t^{(n)}|^{2-2\al}\we 1}\right)\1_{[0,T)}(t)\d t,~r\in [0,s]$$
is a Brownian motion under $\hh R_s\P$ with
\beg{align*}
\hh R_s&=\exp\left\{\int_0^r\left\<\et^{-1}(t)\hh\si^{-1}(t,\hh X_t^{(n),1})\pi_n\left(\ff {\hh X_t^{(n),1}-\hh Y_t^{(n)}} {|\hh X_t^{(n),1}-\hh Y_t^{(n)}|^{2-2\al}\we 1}\right) ,\d W_t\right\>\right.\\
&\qquad \left.-\ff 1 2\int_0^r \left|\et^{-1}(t)\hh\si^{-1}(t,\hh X_t^{(n),1})\pi_n\left(\ff {\hh X_t^{(n),1}-\hh Y_t^{(n)}} {|\hh X_t^{(n),1}-\hh Y_t^{(n)}|^{2-2\al}\we 1}\right) \right|^2\d t \right\}
\end{align*}
and $\hh X_t^{(n),1}$ is a weak solution of \eqref{hhYn} under $\hh R_s\P$. Since $\nn\hh\si\in L^p_q(T)$ with $\ff d p +\ff 2 q<1$, it follows from Theorem \ref{kry-00} that for any $~0\leq s_0<s_1\leq s$
\beg{align*}
\E \hh R_s \left[\int_{r_0}^{r_1} |\cM(\nn\hh\si(t,\cdot))(\hh X_t^{(n),1})|^2\d t\Big|\sF_{s_0}\right]&\leq C\||\cM(\nn\hh\si)|^2\|_{L^{p/2}_{q/2}(s_0,s_1)}\\
&\leq C\| \nn\hh\si \|_{L^{p }_{q }(s_0,s_1)},
\end{align*}
where in the last inequality we use \cite[Lemma 5.4]{ZhangX11}. Then, it follows from  \cite[Lemma 3.5]{XZ} that for any $c\in\R$ 
\beg{align*}
\E \hh R_s\exp\left\{c\int_0^s   |\cM(\nn\hh\si(t,\cdot))(\hh X_t^{(n),1})|^2\d t\right\}<\infty,~s\in [0,T).
\end{align*}
Since $\pi_n$ is bounded, $\E \hh R_s^{-m}<\infty$ for any $m\in\N$ and $s\in [0,T)$. Then the H\"older inequality yields that for any $c$
\beg{align*}
\E \exp\left\{c\int_0^s   |\cM(\nn\hh\si(t,\cdot))(\hh X_t^{(n),1})|^2\d t\right\}<\infty,~s\in [0,T).
\end{align*}
A similar inequality can be established for $\hh X_t^{(n),2}$. Combining these with \eqref{path-ito} and stochastic Gronwall’s inequality, see \cite[Lemma 3.8]{XZ}, the pathwise uniqueness follows. Therefore, the proof is completed.

 \bigskip
 
\noindent\textbf{Acknowledgements}

\medskip

The second author would like to thank Professor Feng-Yu Wang and Dr. Xing Huang for their helpful suggestions. The second author was supported by the National Natural Science Foundation of China (Grant No. 11901604, 11771326).

\end{document}